\documentclass[12pt]{amsart}
\usepackage{amssymb}
\usepackage{cite}
\usepackage{array}
\usepackage{booktabs}
\usepackage{mdwtab}
\usepackage{mathtools}
\usepackage{lmodern}
\usepackage{microtype}
\usepackage[T1]{fontenc}
\usepackage[utf8]{inputenc}
\usepackage{hyphenat}
\usepackage{enumitem}
\usepackage{ifpdf}
\usepackage{braket}
\usepackage{mathabx}
\usepackage{mathrsfs}
\usepackage{color}
\usepackage{minibox} 
\usepackage[pdftex,final]{graphics}
\usepackage[margin=1in]{geometry}
\usepackage[final,bookmarks=true, bookmarksopen=true,%
    bookmarksdepth=3,bookmarksopenlevel=2,%
    colorlinks=true,%
    linkcolor=blue,%
    citecolor=blue,%
    filecolor=blue,%
    menucolor=blue,%
    urlcolor=blue]{hyperref}
  \hypersetup{pdftitle={From Rubber Bands to Rational Maps: A research report}}
  \hypersetup{pdfauthor={Dylan P. Thurston}}
\usepackage[textsize=small]{todonotes}
\usepackage{url}

\usepackage{tikz}
\usetikzlibrary{arrows,calc}
\tikzstyle{every picture}=[> = to]
\tikzset{cdlabel/.style={execute at begin node=$\scriptstyle,execute at end node=$}}
\tikzset{implication/.style={double equal sign distance, -implies}}
\tikzset{biimplication/.style={double equal sign distance, implies-implies}}

\makeatletter
\newcommand\mi@kern[1]{%
  \settowidth\@tempdima{$\mi@obj^{#1}$}
  \kern-\@tempdima
  #1
  \settowidth\@tempdima{$\mi@obj$}
  \kern\@tempdima
}

\newtoks\mi@toksp
\newtoks\mi@toksb
\DeclareRobustCommand{\manyindices}[5]{
  \def\mi@obj{#5}
  \mi@toksp\expandafter{\mi@kern{#2}}
  \mi@toksb\expandafter{\mi@kern{#1}}
  \@mathmeasure4\textstyle{#5_{#1}^{#2}}
  \@mathmeasure6\textstyle{#5_{#3}^{#4}}
  \dimen0-\wd6 \advance\dimen0\wd4
  \@mathmeasure8\textstyle{\hphantom{{}_{#1}^{#2}}#5^{\the\mi@toksp#4}_{\the\mi@toksb#3}}
  \hbox to \dimen0{}{\kern-\dimen0\box8}
}
\makeatother 

\newread\testin

\def\mathcenter#1{\vcenter{\hbox{$#1$}}}
\def\grapha#1{\includegraphics{#1}}
\def\graphb#1{\includegraphics[trim=-1 -1 -1 -1]{#1}}
\def\mfig#1{\mathcenter{\grapha{#1}}}
\def\mfigb#1{\mathcenter{\graphb{#1}}}

\newcommand{\RR}{\mathbb R}
\newcommand{\CC}{\mathbb C}
\newcommand{\DD}{\mathbb D}

\newcommand{\PP}{\mathbb P}

\newcommand{\co}{:}

\newcommand{\abs}[1]{\lvert #1 \rvert}
\newcommand{\norm}[1]{\lVert #1 \rVert}

\newcommand{\bdy}{\partial}

\renewcommand{\Im}{\mathop{\mathrm{Im}}}
\renewcommand{\Re}{\mathop{\mathrm{Re}}}
\DeclareMathOperator*{\esssup}{ess\,sup}
\DeclareMathOperator{\Image}{Im}

\theoremstyle{plain}
\numberwithin{equation}{section}
\newtheorem{proposition}{Proposition}
\numberwithin{proposition}{section}
\newtheorem{lemma}[proposition]{Lemma}
\newtheorem{corollary}[proposition]{Corollary}
\newtheorem{conjecture}[proposition]{Conjecture}

\newtheorem{theorem}{Theorem}
\newtheorem{citethm}[proposition]{Theorem}

\theoremstyle{definition}
\newtheorem{definition}[proposition]{Definition}

\newtheorem{question}[proposition]{Question}
\newtheorem{problem}[proposition]{Problem}

\theoremstyle{remark}
\newtheorem{example}[proposition]{Example}
\newtheorem{remark}[proposition]{Remark}
\newtheorem{warning}[proposition]{Warning}

\DeclareMathOperator{\EL}{EL} 
\DeclareMathOperator{\SF}{SF} 
\DeclareMathOperator{\Fill}{Fill} 

\DeclareMathOperator{\Edges}{Edge}
\DeclareMathOperator{\Lip}{Lip} 
\DeclareMathOperator{\Dir}{Dir} 
\DeclareMathOperator{\Emb}{Emb} 
\DeclareMathOperator{\Teich}{Teich}
\DeclareMathOperator{\Area}{Area}
\renewcommand{\Join}{\mathop{\mathrm{Join}}} 
\newcommand{\ASF}{\overline{\SF}} 
\newcommand{\ALip}{\overline{\Lip}}

\newcommand{\MF}{\mathcal{MF}}
\newcommand{\Curves}{\mathcal{C}}
\newcommand{\Quad}{\mathcal{Q}}
\newcommand{\Wid}{\mathcal{W}}
\newcommand{\Len}{\mathcal{L}}

\newcommand{\powerset}{\mathscr{P}}
\newcommand{\id}{\mathrm{id}}

\newcommand{\wt}[1]{\widetilde{#1}}

\newcommand{\oDD}{\overline{\DD}}

\definecolor{dark-green}{rgb}{0,0.6,0}
\definecolor{dark-red}{rgb}{0.7,0,0}
\definecolor{dark-blue}{rgb}{0,0,0.8}

\graphicspath{{draws/}{figs/}}

\begin{document}
\title[From Rubber Bands to Rational Maps]
  {From Rubber Bands to Rational Maps:\\A Research Report}

\author[Thurston]{Dylan~P.~Thurston}
\thanks{This work was partially supported by NSF grants DMS-1358638 and DMS-1507244.}
\address{Department of Mathematics\\
         Indiana University\\
         831 E. Third St.,
         Bloomington, Indiana 47405\\
         USA}
\email{dpthurst@indiana.edu}

\dedicatory{Dedicated to the memory of William P. Thurston, 1946--2012}

\begin{abstract}
  This research report outlines work, partially joint with
  Jeremy Kahn and Kevin Pilgrim, which gives
  parallel theories of elastic graphs and
  conformal surfaces with boundary. On one hand, this lets us tell
  when one rubber band network is looser than another and, on the
  other hand, tell when one conformal surface embeds in
  another.

  We apply this to give a new characterization of hyperbolic
  critically finite rational maps among branched self-coverings of the
  sphere, by a \emph{positive} criterion: a branched covering is
  equivalent to a hyperbolic rational map if and only if there is an
  elastic graph with a particular ``self-embedding''
  property. This complements the earlier negative criterion of
  W.~Thurston.
\end{abstract}



\maketitle

\setcounter{tocdepth}{1}
  \tableofcontents

\section{Introduction}
\label{sec:intro}

This research report is devoted to explaining a circle of ideas, relating:
\begin{itemize}
\item elastic networks (``rubber bands'') and the corresponding
  Dirichlet energy,
\item extremal length and other ``quadratic'' norms on the space of
  curves in a surface,
\item embeddings between Riemann surfaces,
  conformal or
  quasi-conformal, and
\item post-critically finite rational maps.
\end{itemize}
The original motivation for this project is  the last point, and more
specifically the question of when a topological branched self-covering
of the sphere is equivalent to a rational map. William Thurston first
answered this question more than 30 years ago
\cite{DH93:ThurstonChar}, by giving a ``negative'' characterization: a
combinatorial obstruction that exists exactly when the map is
\emph{not} rational (Theorem~\ref{thm:thurston-obstruction}). In
this paper, we give a ``positive''
characterization: a combinatorial object that exists exactly when the
map \emph{is} rational, for a somewhat restricted
class of maps (Theorem~\ref{thm:detect-rational}).

Having both positive and negative combinatorial
characterizations for the same property automatically gives an
algorithm for testing the property. You search, in parallel, for either an
obstruction or a certificate for the property. One of the searches
will eventually succeed and answer the question. There are other algorithms for
testing whether a branched self-cover is rational: you may take the rational map
itself as a certificate for being
rational
\cite{BBY12:RatlMapDecidable}. However, we
expect our combinatorial certificate to be more practical to search
for (Section~\ref{sec:compute-emb}). By contrast,
W.~Thurston's obstruction theorem is
notoriously hard to apply.

In addition, a positive characterization gives an object to
study associated to the topological branched self-covers of most interest,
namely the rational ones.

Several of the constructions along the way are of independent
interest. For instance, we give a new characterization of when one Riemann
surface conformally embeds in another in a given homotopy class, and
a numerical invariant of such embeddings.

This is a preliminary report on the results. Most proofs are omitted,
although we give some indications. We also indicate
some of the many open problems suggested by this research in
the subsections titled ``Extensions'' (Sects. \ref{sec:extensions},
\ref{sec:extensions-el}, \ref{sec:sf-extensions},
\ref{sec:asf-extensions-questions}, and
\ref{sec:rational-extensions}); all of these are unnecessary for the
main
results.

\subsection{Detecting rational maps}
\label{sec:result}

Recall that a (topological) \emph{branched
  self-cover} of the sphere is a finite set of points~$P$ in the
sphere and a map $f \co (S^2,P) \to (S^2,P)$ that is a covering when
restricted to $S^2 \setminus f^{-1}(P)$. One central question is when
such a map is equivalent to a rational map on $\CC\PP^1$. (See
Definition~\ref{def:branched-self-cover}.) A branched self-cover can
be characterized (up to homotopy) by picking a spine~$\Gamma$ for $S^2
\setminus P$, and drawing its inverse image $\wt\Gamma =
f^{-1}(\Gamma) \subset S^2 \setminus f^{-1}(P)$. There are two natural
homotopy classes
of maps from $\wt\Gamma$ to~$\Gamma$.
(In this paper, a map of graphs is a topological map, not necessarily
taking vertices to vertices.)
\begin{itemize}
\item A covering map commuting with the action of~$f$:
  \[
  \begin{tikzpicture}
    \matrix[row sep=0.8cm,column sep=1cm] {
      \node (Gammai) {$\wt\Gamma$}; &
        \node (Gamma) {$\Gamma$}; \\
      \node (S2i) {$S^2\setminus f^{-1}(P)$}; &
        \node (S2) {$S^2\setminus P$.}; \\
    };
    \draw[->] (Gamma) to node[auto=left,cdlabel] {i} (S2);
    \draw[->] (S2i) to node[auto=left,cdlabel] {f} (S2);
    \draw[->] (Gammai) to node[auto=left,cdlabel] {\pi} (Gamma);
    \draw[->] (Gammai) to (S2i);
  \end{tikzpicture}
  \]
  We denote this covering map $\pi\co \wt\Gamma \to \Gamma$.
\item A map commuting up to homotopy with
  the inclusion in $S^2 \setminus P$:
  \[
  \begin{tikzpicture}
    \matrix[row sep=0.8cm,column sep=1cm] {
      \node (Gammai) {$\wt\Gamma$}; &
        \node (Gamma) {$\Gamma$}; \\
      \node (S2i) {$S^2\setminus f^{-1}(P)$}; &
        \node (S2) {$S^2\setminus P$.}; \\
    };
    \draw[->] (Gamma) to node[auto=left,cdlabel] {i} (S2);
    \draw[right hook->] (S2i) to node[auto=left,cdlabel] {\text{id}} (S2);
    \draw[->] (Gammai) to node[auto=left,cdlabel] {\phi} (Gamma);
    \draw[->] (Gammai) to (S2i);
  \end{tikzpicture}
  \]
  We denote this map~$\phi$. It is unique up to homotopy since
  $\Gamma$ is a spine for $S^2\setminus P$.
\end{itemize}

\begin{figure*}
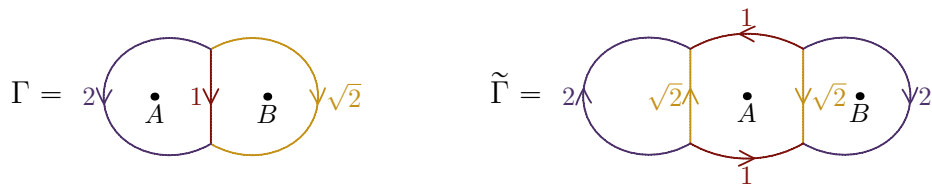

  \centerline{
    $\Gamma    = \,\,\,\,\mfigb{graphs-11}$\hfil
    $\wt\Gamma = \,\,\,\,\mfigb{graphs-12}$}
  \caption{A spine and its inverse image defining a branched
    self-cover of the sphere. The marked set~$P$ consists of the two
    points shown, and a third point~$C$ at infinity.  The numbers on
    the edges give a measure which proves this branched self-cover is
    equivalent to a rational map.}
  \label{fig:spine-example}
\end{figure*}
\begin{example}\label{examp:3points-map}
  Figure~\ref{fig:spine-example} shows a simple example of such spines
  when $P$ has 3 points, $A$, $B$, and~$C$ (at~$\infty$). The covering
  map~$\pi$ preserves the colors and arrows on the edges.
  The branched self-cover $f\co (S^2,P)\righttoleftarrow$ is an
  extension of $\pi$
  to the
  whole sphere. This extension is unique up to homotopy relative to
  $\{A,B,C\}$. The map $f$
  permutes the marked points by
  \[
  \begin{tikzpicture}
    \node (A) at (150:0.8cm) {$A$};
    \node (B) at (30:0.8cm) {$B$};
    \node (C) at (-90:0.8cm) {$C$};
    \draw[|->,bend left=20] (A) to node[above,cdlabel] {(2)} (B);
    \draw[|->,bend left=20] (B) to (C);
    \draw[|->,bend left=20] (C) to node[below left=-2pt,cdlabel] {(2)} (A);
  \end{tikzpicture}
  \]
  (For instance, $A$ is inside a crimson-goldenrod region in~$\wt\Gamma$, so
  must map to $B$, which is inside a crimson-goldenrod region in~$\Gamma$.)

  The map $\phi$, on the other hand, is the projection of
  $\wt\Gamma$ onto $\Gamma$ considered as a spine; for instance, it
  might map the right-hand purple edge of $\wt\Gamma$ to the right-hand
  goldenrod edge of $\Gamma$, whereas $\pi$ preserves the colors.
\end{example}

For our characterization of rational maps, we also need an
\emph{elastic structure} on~$\Gamma$, by which we mean
a measure~$\alpha$ on~$\Gamma$ absolutely continuous with respect to
Lebesgue measure.%
\footnote{We can also interpret $\alpha$ as
  a metric. But we prefer to distinguish this measure from other
  metrics that come later.}
We can pull back $\alpha$ by~$\pi$ to get a measure
$\pi^* \alpha = \wt\alpha$ on~$\wt\Gamma$.
For $\psi\co \wt\Gamma \to \Gamma$ a Lipschitz map
of graphs, define the \emph{embedding energy} by
  \begin{equation}
    \label{eq:embedding-1}
    \Emb(\psi) \coloneqq
      \esssup_{y \in \Gamma} \sum_{\substack{x \in \wt\Gamma\\\psi(x)=y}} \abs{\psi'(x)}.
  \end{equation}
The derivative is taken with respect to the two measures $\wt\alpha$
and~$\alpha$. The \emph{essential supremum} $\esssup$ ignores points
of measure zero.  In particular, we may ignore
vertices of~$\wt\Gamma$ or points of $\wt\Gamma$ that map to vertices
of~$\Gamma$.

In practice the embedding energy is optimized when $\psi$ is
piecewise linear, and the reader may restrict to that case.

For one motivation for
Equation~\eqref{eq:embedding-1}, see
Proposition~\ref{prop:sf-graph-surface}: the embedding energy
characterizes conformal embeddings of thickened versions of the
graph. Another one is in Theorem~\ref{thm:sf-emb-energy}: it
characterizes when Dirichlet energy is reduced under the map,
i.e., when one rubber band network is looser than another. (See
Remark~\ref{rem:looser}.)

\begin{example}\label{examp:3points-emb}
  To return to Example~\ref{examp:3points-map}, consider the measure $\alpha$
  on~$\Gamma$ from Figure~\ref{fig:spine-example} and the concrete
  map~$\psi\co \wt\Gamma \to \Gamma$ that maps
  \begin{itemize}
  \item the right purple edge of $\wt\Gamma$ to the goldenrod edge on the
    right of $\Gamma$,
  \item the right goldenrod edge of $\wt\Gamma$ to the crimson edge in the
    middle of $\Gamma$, and
  \item the remaining four edges (two crimson, one goldenrod, and one purple) of
    $\wt\Gamma$ map to
    the purple edge on the right of $\Gamma$, with the crimson edges mapping to
    segments of length $1/\sqrt{2}$.
  \end{itemize}
  Make $\psi$ linear on each segment described above.
  Then a short computation shows that $\Emb(\psi) = 1/\sqrt{2}$.
\end{example}

A major theorem of this paper is that rational maps can be
characterized by maps~$\psi$ with
embedding energy less than~$1$.
We say that two branched self-covers
$f_1\co (S^2,P_1)\righttoleftarrow$ and
$f_2 \co (S^2,P_2)\righttoleftarrow$ are \emph{equivalent} if they can
be connected by a homotopy of $f_1$ relative to $P_1$ and conjugacy
taking $P_1$ to $P_2$.

\begin{theorem}\label{thm:detect-rational}
  Let $f \co (S^2,P) \righttoleftarrow$ be a branched self-cover of
  the sphere. Suppose that there is a branch point in each cycle
  in~$P$. Then  $f$
  is equivalent to a rational map iff there is
  an elastic graph spine~$(\Gamma,\alpha)$ for $S^2 \setminus P$,
  an integer $n > 0$, and map
  $\psi\in[\phi_n]$ so that $\Emb(\psi) < 1$.
\end{theorem}
Here $\phi_n$ is obtained from $\phi$ by iteration, and is in the
homotopy class of the projection
$f^{-n}(\Gamma)\to\Gamma$.
Loosely speaking, Theorem~\ref{thm:detect-rational} says that the
self-cover is rational iff there is a
\emph{self-embedded spine} for $S^2\setminus P$.

\begin{remark}\label{rem:detect-relax}
  It is likely the condition on Theorem~\ref{thm:detect-rational} can
  be relaxed to assume merely that $f$ has at least one branch point
  in one cycle in~$P$, i.e., if $f$ is rational, its Julia set is not
  the whole sphere. See Section~\ref{sec:rational-extensions}.
\end{remark}

\begin{example}
  The explicit measure and map in Example~\ref{examp:3points-emb} show
  that the given branched self-cover is equivalent to a rational
  map. In fact, it is equivalent to the rational map
  \[
  f(z) = \frac{1}{1-z^2},
  \]
  with $A = 0$, $B = 1$, and $C = \infty$.%
\footnote{Every branched
  self-cover with only $3$ post-critical points is equivalent to a
  rational map, so this example was easy to do by other means.}
\end{example}

There was a previous characterization of rational maps by W.~Thurston,
recalled as Theorem~\ref{thm:thurston-obstruction} below.
This is analogous to the two ways to characterize pseudo-Anosov
surface automorphisms, which form a natural class of geometric
elements of the mapping class group of a surface. Geometrically,
pseudo-Anosov diffeomorphisms are those whose mapping torus is a
hyperbolic 3-manifold. Combinatorially, there are two
criteria:
  \begin{itemize}
  \item \emph{Negative}: $f$ is pseudo-Anosov iff it is not
    periodic ($f^{\circ k}$ is not the identity for any $k > 0$) or
    reducible (there is no invariant system of multi-curves for~$f$).
  \item \emph{Positive}: $f$ is pseudo-Anosov iff there is a
    measured train track~$T$ and a splitting sequence from $T$ to a
    train-track~$T'$ with $T' = \lambda f(T)$ for some constant $0 <
    \lambda < 1$ \cite{PP87:CharPAnosov}.
  \end{itemize}
The positive criterion gives some extra information: the number
$\lambda$ is an invariant of~$f$, with dynamical interpretations. (For
instance, $\lambda$ controls the growth rate of intersection numbers.)

Analogously, we can say that a branched self-cover~$f$ of $(S^2,P)$ is
geometric if it is equivalent to a rational map, which also
have associated 3-dimensional hyperbolic laminations
\cite{LM97:LaminationsHolomorphic}. There are combinatorial criteria
for $f$ to be rational:
\begin{itemize}
\item \emph{Negative}: $f$ is rational iff there is no
  obstruction, as in W.~Thurston's
  Theorem~\ref{thm:thurston-obstruction}. In loose terms, the
  obstruction is a back-expanding annular system: a collection of
  annuli that get ``wider'' under backwards iteration.
\item \emph{Positive}: Under some additional assumptions, $f$ is
  rational iff there is a
  metric spine for~$S^2\setminus P$ satisfying a back-contracting
  condition. This is Theorem~\ref{thm:detect-rational}.
\end{itemize}
As in the case of surface automorphisms, the two theorems are in a
sense ``dual'' to each other: it is easy to see that a branched-self
cover cannot simultaneously have a back-expanding annular system and a
back-contracting spine. (See Equation~\eqref{eq:obstruction-cd}.) Also
as in the surface automorphism case, the
positive criterion gives us a new object to study, namely the constant
$\ASF[\phi]$ of
Section~\ref{sec:asymptotic-SF}.

Compared to the situation for surface automorphisms,
Theorem~\ref{thm:detect-rational} has the following caveats:
\begin{itemize}
\item It only works when there is a branch point in each cycle
  in~$P$. (But see Remark~\ref{rem:detect-relax}.)
\item We may need to pass to an iterate to get $\Emb(\phi) <
  1$. Furthermore, it is easy to see in examples that the embedding
  energy can decrease in powers. (In the notation of
  Section~\ref{sec:asymptotic-SF}, in general $\SF[\phi_n] <
  \SF[\phi]^n$ and $\ASF[\phi] < \SF[\phi]$. See
  Example~\ref{examp:asf-bad}.)
\end{itemize}

There is a further caveat for both surface automorphisms and branched
self-covers:
\begin{itemize}
\item For the positive criterion, the train track or graph constructed
  is not canonical: there are many different choices that work for the
  criterion.
\end{itemize}
 By contrast, the negative criteria can be made
  canonical. (Pilgrim \cite{Pilgrim01:CanonicalObstruction} proved
  this for branched self-covers.)
On the other hand, in the surface case, Agol
\cite{Agol11:pAtriangulation} and Hamenstädt
\cite{Hamenstaedt09:GeomMCG1} give a canonical object related to the measured
train track.

\subsection{Elastic graphs}
\label{sec:overview}

The embedding energy of Equation~\eqref{eq:embedding-1} looks a little
mysterious; it looks a little like the Lipschitz stretch factor, but
the sum over inverse images looks unusual. To explain where it comes
from, we now turn to a ``conformal'' theory of graphs
parallel to the conformal theory of Riemann surfaces. The central
object is an \emph{elastic graph}~$(\Gamma,\alpha)$, which you should
think of as a network of rubber bands; formally, it is a graph with a
measure on each edge, representing the elasticity of the edge. (See
Section~\ref{sec:graphs}.)

There are additional structures we can put
on the graph.
\begin{itemize}
\item On one hand, we can consider \emph{curves}~$C$ on the graph,
  maps
  from a 1-manifold into~$\Gamma$.
\item On the other hand, we can consider maps from~$\Gamma$ to a
  \emph{length graph}~$K$, a graph
  with fixed lengths of edges (like a network of pipes).
\end{itemize}
Maps between these objects have naturally associated energies, as
summarized in the following diagram.
\begin{equation}
\mathcenter{\begin{tikzpicture}[node distance=3cm,
  every text node part/.style={align=center},align=center]
  \node (curve) {Curve\\ $C$};
  \node (elastic) [right of=curve] {Elastic\\graph $\Gamma$};
  \node (length) [right of=elastic] {Length\\graph $K$};
  \draw[->] (curve) to node[above]{$\EL$} (elastic);
  \draw[->] (elastic) to node[above]{$\Dir$} (length);
  \draw[->, bend left=35] (curve) to node[above]{$\ell^2$} (length);
  \draw[->, loop below, min distance=1cm, out=-70, in=-110]
    (elastic) to node[below]{$\Emb$} (elastic);
  \draw[->, loop below, min distance=1cm, out=-70, in=-110]
    (length) to node[below]{$\Lip^2$} (length);
\end{tikzpicture}}
\label{eq:energies-graph}
\end{equation}
The labels on the arrows indicate the type of energy on a map of this
type, as follows.
\begin{itemize}
\item For a map~$f$ from an elastic graph $(\Gamma,\alpha)$ to a
  length graph $(K,\ell)$,
  there is the \emph{Dirichlet} or \emph{rubber-band energy}
  (Section~\ref{sec:harmonic}) familiar from physics:
  \begin{equation}
    \label{eq:dirichlet-1}
    \Dir(f) \coloneqq \int_{x \in \Gamma} \abs{f'(x)}^2\,dx,
  \end{equation}
  where $f'$ measures the derivative with respect to the natural
  metrics. If $f$ minimizes this energy within its homotopy class, it
  is said to be \emph{harmonic}.
\item For a curve~$C$ in an elastic graph $(\Gamma,\alpha)$, we have a
  version of
  \emph{extremal length} (Section~\ref{sec:extremal-length}):
  \begin{equation}
    \label{eq:EL-1}
    \EL[C] \coloneqq \sum_{e \in \Edges(\Gamma)} n_C(e)^2\cdot \alpha(e),
  \end{equation}
  where $n_C(e)$ is the number of times $C$ runs over the edge~$e$
  (without backtracking).
\item For a curve~$C$ in a length graph $(K,\ell)$, we have the
  usual \emph{length}, which in our notation is
  \begin{equation}
    \label{eq:length-1}
    \ell[C] \coloneqq \sum_{e \in \Edges(K)} n_C(e) \cdot \ell(e).
  \end{equation}
  To match the other quantities, we actually use the square of
  the length as our energy.
\item For a map~$\phi$ between length graphs $(K_1,\ell_1)$ and
  $(K_2,\ell_2)$, there is the \emph{Lipschitz constant}
  (Section~\ref{sec:lipschitz-energy}):
  \begin{equation}
    \label{eq:lipschitz-1}
    \Lip(\phi) \coloneqq \esssup_{x \in K_1} \abs{\phi'(x)}.
  \end{equation}
  Again, we consider the square of the Lipschitz energy.
\item Finally, for a map~$\phi$ between elastic
  graphs $(\Gamma_1,\alpha_1)$ and $(\Gamma_2, \alpha_2)$, we have
  \emph{embedding energy} as used in Theorem~\ref{thm:detect-rational}
  (Section~\ref{sec:stretch-factors}):
  \begin{equation}
    \label{eq:embedding-2}
    \Emb(\phi) \coloneqq
      \esssup_{y \in \Gamma_2} \sum_{x \in \phi^{-1}(y)} \abs{\phi'(x)}.
  \end{equation}
  This energy appears to be new, although it is related to Jeremy
  Kahn's notion of domination of weighted arc diagrams
  (Section~\ref{sec:dynam-teich}).
\end{itemize}
\begin{remark}
We could make the diagram more symmetric by using \emph{width graphs}
instead of curves (see Section~\ref{sec:graphs}), and adding a norm on
maps between width graphs.
\end{remark}

These energies are \emph{sub-multiplicative}, in the sense that
composing two maps can only decrease the product of the
energies: if $f$ and $g$ are two composable maps of the above types,
then
\begin{equation}
  \label{eq:norm-submultiplicative}
  \norm{f \circ g} \le \norm{f}\cdot\norm{g},
\end{equation}
where $\norm{\cdot}$ is the appropriate energy from the above list.
(This inequality is the reason we squared some of the energies.) For
instance, if we fix elastic graphs $\Gamma_1$ and~$\Gamma_2$, a
length graph~$K$, and maps $\phi \co \Gamma_1 \to \Gamma_2$ and $f \co
\Gamma_2 \to K$, then
\begin{equation}\label{eq:emb-inequality-weak}
  \Dir(f \circ \phi) \le \Dir(f) \cdot \Emb(\phi)
\end{equation}
What is more,
these inequalities are all \emph{tight}, in the sense that if we fix the
domain, range, and homotopy type of $f$, then we can find a
sequence of functions~$g_i$ (including a choice of domain)
that approach
equality in Equation~\eqref{eq:norm-submultiplicative}.  Likewise
if we fix~$g$ and vary~$f$, we can find a sequence of functions~$f_i$
approaching equality in Equation~\eqref{eq:norm-submultiplicative}.

For instance, if $\phi \co \Gamma_1 \to \Gamma_2$ is a map of elastic
graphs, we can strengthen Equation~\eqref{eq:emb-inequality-weak} to
\begin{equation}\label{eq:emb-equality}
  \Emb[\phi] =
    \sup_{\substack{\text{$K$ length graph}\\f: \Gamma_2 \to K}}
      \frac{\Dir[f \circ \phi]}{\Dir[f]},
\end{equation}
where $\Emb[\phi]$ and $\Dir[f]$ are the minimums over the respective
homotopy classes (Theorem~\ref{thm:sf-emb-energy}). Since Dirichlet
energy can be interpreted as the
elastic energy of a stretched rubber band network, $\Emb[\phi] < 1$ 
can therefore also be interpreted as saying that $\Gamma_1$ is
``looser'' than $\Gamma_2$, however the two rubber band networks
are stretched. See Remark~\ref{rem:looser}.

Other examples are in
Propositions~\ref{prop:dir-el-dual} and~\ref{prop:el-dir-dual}.
This gives a kind of duality between curves in an elastic
graph~$\Gamma$ and maps from~$\Gamma$ to
length graphs. If we think of curves as living in a ``vector space''
and maps to length graphs in its ``dual'', then the embedding energy can be
interpreted as an ``operator norm''.

This theory of conformal graphs is largely parallel to the theory
of conformal (Riemann) surfaces with boundary, where we again have a number of
energies:
\begin{equation}
\mathcenter{\begin{tikzpicture}[node distance=3cm,
  every text node part/.style={align=center},align=center]
  \node (curve) {Curve\\ $C$};
  \node (surf) [right of=curve] {Riemann\\surface $\Sigma$};
  \node (length) [right of=elastic] {Length\\graph $K$};
  \draw[->] (curve) to node[above]{$\EL$} (surf);
  \draw[->] (surf) to node[above]{$\Dir$} (length);
  \draw[->, bend left=35] (curve) to node[above]{$\ell^2$} (length);
  \draw[->, loop below, min distance=1cm, out=-70, in=-110]
    (surf) to node[below]{$\SF$} (surf);
  \draw[->, loop below, min distance=1cm, out=-70, in=-110]
    (length) to node[below]{$\Lip^2$} (length);
\end{tikzpicture}}\label{eq:energies-surface}
\end{equation}
Again, each arrow is marked by the appropriate energy for measuring
that type of map.
$\Dir$ and $\EL$ are again Dirichlet energy and extremal length,
but on surfaces rather than graphs. $\SF$ is new; it is the
\emph{stretch factor}
of a homotopy class of a topological embedding $[\phi]\co \Sigma_1 \to
\Sigma_2$ between Riemann surfaces. 
In general, we do not know a direct expression analogous to
Equation~\eqref{eq:embedding-1}, so $\SF$ is defined to be the minimal
ratio of extremal lengths
(Definition~\ref{def:sf-surface}, analogous to
Equation~\eqref{eq:emb-equality}). When
there is no conformal embedding of
$\Sigma_1$ in~$\Sigma_2$ in the given homotopy class,
there is a direct expression: $\SF[\phi]$ is given by the
minimal quasi-conformal constant in the homotopy class
(Theorem~\ref{thm:SF-Teich}).
There is an analogue of $\SF$ for maps between elastic graphs, and in
that context $\SF[\phi] = \Emb[\phi]$ (Theorem~\ref{thm:sf-emb-energy}).

\begin{remark}
To prove Theorem~\ref{thm:detect-rational}, we  actually do
not need to consider length graphs or the Dirichlet energy at
all. (The proofs go through extremal length instead.) However, they
illuminate the overall structure. In
particular, it is not clear why one would consider elastic graphs
without the rubber-band motivation.
\end{remark}

\begin{remark}
  The appearance of length squared in \eqref{eq:energies-graph} and
  \eqref{eq:energies-surface} is easy to justify on the grounds of
  units. Extremal length itself behaves like the square of a length,
  in the sense that if we take $k$ parallel copies of a curve, the
  extremal length multiplies by $k^2$. Likewise, if the lengths on the
  target of a harmonic map are multiplied by $k$, the harmonic map
  remains harmonic while the Dirichlet energy
  is multiplied by $k^2$.
\end{remark}

\subsection{History and prior work}
\label{sec:history}

Although Equation~\eqref{eq:embedding-1} appears to be
new, Jeremy Kahn's notion of \emph{domination} of weighted arc
diagrams \cite{Kahn06:BoundsI} is essentially equivalent. See
Section~\ref{sec:dynam-teich}.

There has also been substantial work in
the related setting of resistor networks rather than spring networks
\cite[inter alia]{BSST40:DissectionRects,Duffin62:ELNetwork,CIM98:CircNetworks}.
See Remark~\ref{rem:not-electrical}.

Theorem~\ref{thm:detect-rational} is closely related to
Theorem~\ref{thm:rational-surfaces-embed}, which characterizes when a
rational map exists in terms of conformal embeddings of
surfaces. Theorem~\ref{thm:rational-surfaces-embed} has
been a folk theorem in the community for some time.

For polynomials, Theorem~\ref{thm:detect-rational} reduces to a
previously-known characterization in terms of expansion on the Hubbard
tree; see Theorem~\ref{thm:realize-hubbard} in
Section~\ref{sec:polynomials}.

\subsection{Organization}
\label{sec:organization}

After a section giving some examples of how to apply
Theorem~\ref{thm:detect-rational}, this paper is organized by
topics moving up a dynamical hierarchy. For each topic we give first
the
conformal surface notions and then the graph notions.
\begin{itemize}
\item We start with notions depending only on a single conformal
  surface or graph: \emph{Dirichlet} (or \emph{rubber band})
  \emph{energy} of (harmonic) maps (Section~\ref{sec:harmonic}) and
  \emph{extremal length} of curves (Section~\ref{sec:extremal-length}).
\item Next come notions depending on a map between surfaces or
  graphs. This is the \emph{stretch factor} or \emph{embedding
    energy}, which generalizes the Teichm\"uller distance and
  characterizes conformal embeddings of surfaces
  (Section~\ref{sec:stretch-factors}).
\item Next is the dynamical theory of iterated maps
  (Section~\ref{sec:dynamics}). Here we find another number, the
  \emph{asymptotic stretch factor}, which characterizes rational maps
  (Section~\ref{sec:rational-maps}).
\end{itemize}
The paper is organized by a logical hierarchy, rather than what is
necessarily pedagogically best; the reader is encouraged to skip
around.

\subsection{Acknowledgements}
\label{sec:acknowledgements}

I would like to thank
Matt Bainbridge,
Steven Gortler,
Richard Kenyon,
Sarah Koch,
Tan Lei,
Dan Margalit, and
Giulio Tiozzo
for many helpful conversations.
I would like to especially thank Maxime Fortier Bourque, who pointed
me towards Ioffe's theorem
\cite{Ioffe75:QCImbedding} and had numerous other insights.

This project grew out of extensive conversations with Kevin Pilgrim,
who helped shape my understanding of the subject in many ways. Many of
the arguments were developed jointly with him. Notably, he
communicated Theorem~\ref{thm:rational-surfaces-embed} to me, and
Theorem~\ref{thm:el-embedding} is joint work with him.
Theorem~\ref{thm:SF-cover-bound} is joint work with Jeremy
Kahn, who also contributed substantially throughout.

Above all, I would like to thank William Thurston for
introducing me to the subject and insisting on
understanding deeply.


\section{Examples}
\label{sec:examples}

Here, we give some more substantial examples of
Theorem~\ref{thm:detect-rational}.

\subsection{Polynomials: The rabbit and the basilica}
\label{sec:poly-examples}

Theorem~\ref{thm:detect-rational} is not very interesting for
polynomials, as \emph{every} topological polynomial with a
branch-point in each cycle is equivalent to a polynomial. The
extension of Theorem~\ref{thm:detect-rational} to the general
topological polynomial case is somewhat more interesting, but is
equivalent to known results on expansion on the Hubbard tree; see
Section~\ref{sec:polynomials}. Nevertheless, we will look at some
examples, both to see what the stretch factors are and to use them for
matings.

\begin{example}\label{examp:rabbit-1}
We first look at the ``rabbit'' polynomial, the post-critically finite
polynomial $f_1(z) = z^2+c$ with $c \approx -0.1226 + 0.7449i$. The
critical point moves in a 3-cycle
\[
\mathcenter{\begin{tikzpicture}[node distance=1.5cm]
  \node (A) at (-90:0.9cm) {$0$};
  \node (B) at (30:0.9cm) {$c$};
  \node (C) at (150:0.9cm) {$c^2 + c$};
  \draw[|->,bend right=15] (A) to node[right,cdlabel,pos=0.4]{(2)} (B);
  \draw[|->,bend right=15] (B) to (C);
  \draw[|->,bend right=15] (C) to (A);
\end{tikzpicture}}.
\]
The optimal elastic graph $\Gamma_1$ and its cover $\wt\Gamma_1$ are
\[
\mfigb{matings-51} \longrightarrow \;\;\mfigb{matings-50},
\]
with $\Emb[\wt\Gamma_1 \to \Gamma_1] = 2^{-1/3}$, which is less than one,
as expected. (There is another marked point at infinity, not shown.)
\end{example}

\begin{example}\label{examp:rabbit-2}
Another graph $\Gamma_2$ that works to prove that the rabbit polynomial is
realizable is
\begin{equation}
  \label{eq:rabbit-2}
   \mfigb{matings-2} \;\;\longrightarrow \;\;\mfigb{matings-1}.
\end{equation}
The black edges have the indicated lengths, which come from looking at
the external rays landing at the $\alpha$ fixed point of $f_1$. Give
the colored edges an equal and
long elastic length (say, 100). There is a natural map $\phi_0 \co
\wt\Gamma_2 \to \Gamma_2$ as follows. 
\begin{itemize}
\item The outside circle is mapped to the outside circle, with derivative
  $\abs{\phi_0'} = 1/2$.
\item The colored segments on the lower right of $\wt\Gamma_2$ is squashed
  out to the lower-right boundary of $\Gamma_2$, with derivative
  $\abs{\phi_0'}$ on order of $1/100$. Thus $\Fill_{\phi_0} \approx 1/2$ on
  the corresponding portion of $\Gamma_2$.
  ($\Fill_\phi$ is defined in Equation~\eqref{eq:fill-def}, and is
  the quantity maximized in the definition of $\Emb(\phi)$.)
\item The colored segments in the upper left of $\wt\Gamma_2$ are
  mapped to the
  colored segments in $\Gamma_2$, with $\abs{\phi_0'} = 1$.
\end{itemize}
Because of the last point, this map has $\Emb[\phi_0] = 1$, which is
not good enough. Let $\phi_1$
be the result of ``pulling in'' very slightly the image of the ends of the
upper-left colored segments of $\wt\Gamma_2$, so they map to the
interior of the colored segments of $\Gamma_2$. This decreases the
derivative on the colored segments to less than one, while only
increasing the derivative on the outside circle slightly; thus
$\Emb[\phi_1]$ is very slightly less than~$1$, as desired.
\end{example}

\begin{example}\label{examp:basilica}
We can perform a similar trick for other polynomials. For instance, the
basilica polynomial $f_3(z) = z^2-1$ has a spine $\Gamma_3$, with cover
given by
\begin{equation}
  \label{eq:basilica}
  \mfigb{matings-11} \longrightarrow \mfigb{matings-10}.
\end{equation}
The same argument as above (giving the purple edge a long length,
and pushing the right purple edge out to the boundary) shows that
there is a map $\phi \co \wt\Gamma_3 \to \Gamma_3$ with $\Emb[\phi] <
1$. (In this case the optimal stretch factor is
$2^{-1/2}$ and is not realized with a graph with this topology).
\end{example}

\subsection{Matings}
\label{sec:matings}

We can use the techniques of Section~\ref{sec:poly-examples} to show
that some matings of polynomials are geometrically
realizable.

\begin{example}\label{examp:mate-rabbit-basilica}
We can glue together the figures in
Equations~\eqref{eq:rabbit-2} and~\eqref{eq:basilica}:
\begin{equation}
  \label{eq:mate-rabbit-basilica}
  \mfigb{matings-21} \longrightarrow \mfigb{matings-20}.
\end{equation}
This gives a graph spine $\Gamma_4$ and cover $\wt\Gamma_4$
representing the formal mating of the rabbit and the basilica.
We can find a map $\phi \co \wt\Gamma_4 \to \Gamma_4$ with embedding
energy less than $1$:
\begin{itemize}
\item Assign the black mating circle in $\Gamma_4$ total length~$1$,
  divided according to the angles of the external rays.
\item Give all colored edges an equal and large length.
\item Pull this metric back to $\wt\Gamma_4$, and map $\wt\Gamma_4$ to
  $\Gamma_4$ by pushing colored arcs out to the black circle as in
  Examples~\ref{examp:rabbit-2} and \ref{examp:basilica}.
\item Pull the map in slightly where colored vertices meet the black
  circle.
\end{itemize}
The result has $\Fill_\phi(y)$ slightly larger than $1/2$ when $y$ is
on the black circle, and $\Fill_\phi(y)$ slightly less than one on the
colored edges.
\end{example}

Naturally, the technique of Example~\ref{examp:mate-rabbit-basilica}
cannot always work, as sometimes the mating is not
geometrically realizable.

\begin{example}\label{examp:mate-basilica-basilica}
If we try to mate a basilica with a
basilica, we get these graphs:
\begin{equation}
  \label{eq:mate-basilica-basilica}
  \mfigb{matings-31} \longrightarrow \mfigb{matings-30}.
\end{equation}
If we try to use the same technique as before, it doesn't work, as
there are two points on the black circle that we attempt to pull in
two different directions. Indeed, the left green-purple circle is
mapped to a green-purple circle, so must have derivative at least
$1$: it is an obstruction to the mating. (In this case, it is a Levy
cycle).
\end{example}

\subsection{Slit maps}\label{sec:slit-examples}
Given a branched self-cover $f\co (S^2,P)\righttoleftarrow$ and an
arc~$A$ with endpoints in~$P$, there
is a \emph{blowing up} construction which produces a map $f_A$ that agrees
with $f$ outside of a neighborhood of~$A$ and maps that neighborhood
surjectively on to~$S^2$.
Pilgrim and Tan Lei showed that, if $f$ is a rational map, these blow
ups frequently are as well \cite{PT98:Combining}. We will restrict
attention to cases where the initial map $f$ is the identity, in which
case the theorem becomes the following.

\begin{citethm}\label{thm:slit}
  Let $P \subset S^2$ be a finite graph, and let $G \subset S^2$ be a
  finite embedded graph with endpoints on~$P$. Then $\id_G$
  is a rational map iff $G$ is connected.
\end{citethm}

We can give a new proof in the harder direction, when $G$ is
connected. We start with a simple example. From the connected planar graph
\[
G = \mfigb{graphs-40},
\]
take as spine~$\Gamma_4$ the spherical dual to~$G$. Then $\wt\Gamma_4$ is
obtained by taking the connect sum of~$\Gamma_4$ with four extra
copies $\Gamma_4^i$
of~$\Gamma_4$, one for each edge:
\[
\mfigb{graphs-42}\longrightarrow\mfigb{graphs-41}.
\]
For any metric on $\Gamma_4$, there is a natural map $\phi \co
\wt\Gamma_4 \to \Gamma_4$ that maps most of each copy $\Gamma_4^i$ to a point
in the center of the corresponding edge of $\Gamma_0$. This map
$\phi$ has
derivative equal to $1/2$ or $0$ everywhere. We have $\Emb(\phi) =
1/2$, and consideration of the red edge shows that this is optimal.

This example generalizes immediately to a general connected graph~$G$,
except that $\Emb[\phi]$ will only be $1/2$ when $G$ has a univalent
vertex; otherwise, $\Emb[\phi]$ will be strictly smaller. We have thus
proved Theorem~\ref{thm:slit}, with some additional information about
the stretch factor.

\subsection{Behavior under iteration}
\label{sec:more}

We can iterate Example~\ref{examp:3points-map}:
\[
\cdots \rightarrow \mfigb{graphs-3} \rightarrow \mfigb{graphs-2}
\rightarrow \mfigb{graphs-1}
\]
(See Section~\ref{sec:iterating-covers} for details on iteration.) In
this example, the embedding energy behaves well, in the
sense that the $k$'th iterate $\phi_k$ is
optimal for embedding energy (using the edge lengths and concrete
initial map from
Example~\ref{examp:3points-emb}):
\[
\Emb[\phi_k] = \Emb(\phi_k) = \bigl(\sqrt{2}\bigr)^{-k}.
\]
Such good behavior is not generally the case.

\begin{example}\label{examp:asf-bad}
  The map
  \[
  f(z) = \frac{1+z^2}{1-z^2}
  \]
  is represented combinatorially by a graph $\Gamma_5$ and maps
  $\pi,\phi\co \wt\Gamma_5 \rightrightarrows \Gamma_5$:
  \begin{equation}\label{eq:asf-bad}
    \mfigb{graphs-5}\longrightarrow
    \mfigb{graphs-1}.
  \end{equation}
  A case analysis
  shows that for any metric on $\Gamma$ and any map $\psi \in [\phi]$ with
  $\Emb(\psi) < 1$, there is some power $1 \le k \le 4$ so that the
  $k$-fold iterate $\psi_k$ has local back-tracking. This implies that
  \[
  \Emb[\psi_4] < \Emb[\psi]^4,
  \]
  regardless of the initial metric.
  Similar facts hold for any of the other three graphs
  homotopy equivalent to~$\Gamma_5$.
  Thus, regardless of the choice of initial spine, we have $\ASF[\phi]
  < \SF[\phi] = \Emb[\phi]$. (See
  Section~\ref{sec:asymptotic-SF} for definitions.)

  Figure~\ref{fig:computations} gives a sample of experimental data for
  this map. An ad hoc argument shows that for this map, $\ASF[\phi] =
  2^{-1/3}$. The corresponding
  line is shown dashed in the figure.
  \begin{figure}
    \centering
    \includegraphics[scale=0.8]{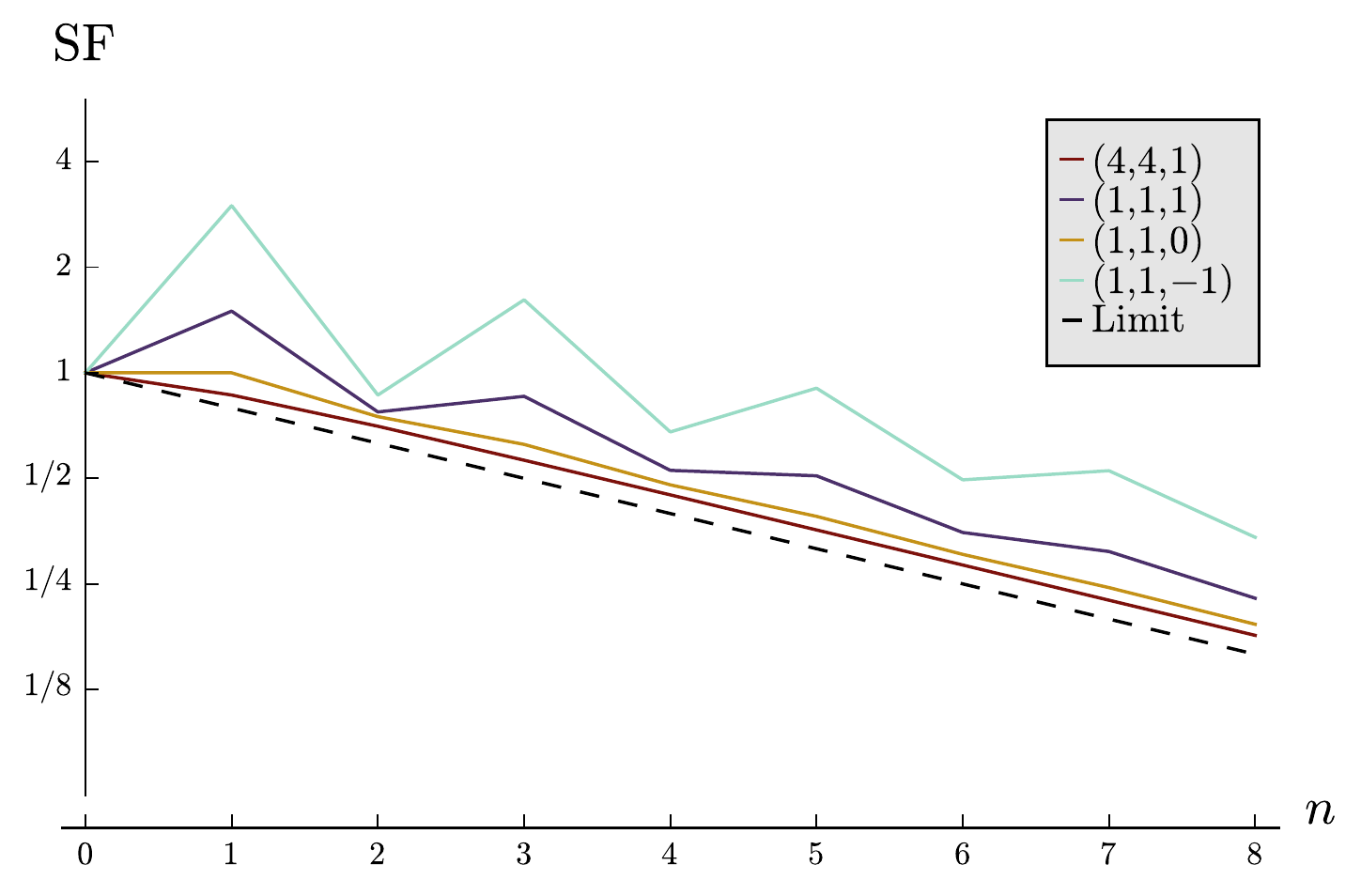}%
    \caption{Logarithmic plot of $\SF[\phi_n]$ for $\phi$ from
      Example~\ref{examp:asf-bad}, as $n$ varies. The parameters
      $(\alpha,\beta,\gamma)$ are the elastic lengths of the left, right, and
      middle edges in Equation~\eqref{eq:asf-bad}, respectively. When
      $\gamma=0$, the graph $\Gamma_5$ degenerates to a rose
      graph. The line for $(\alpha,\beta,\gamma)= (1,1,-1)$
      corresponds instead to a barbell graph
      homotopy equivalent to $\Gamma_5$. The dashed line shows $\ASF[\phi]^n$.}
    \label{fig:computations}
  \end{figure}
\end{example}


\section{Setting}
\label{sec:setting}

\subsection{Surfaces}
\label{sec:surfaces}

We work with compact, oriented surfaces~$\Sigma$ with boundary.
It is sometimes convenient to think about the \emph{double} $D\Sigma
= \Sigma \cup_\partial \overline{\Sigma}$, which has no boundary.

\begin{definition}
A \emph{curve} on~$\Sigma$ is an immersion of a 1-manifold with
boundary in~$\Sigma$, with boundary mapped to the boundary. The
1-manifold need not be connected; if it is, the curve is said to be
\emph{connected}. Curves are
considered up to homotopy within the space of all maps taking the
boundary to the boundary (not
necessarily immersions). A curve
is \emph{simple} if it is embedded (has no crossings). An \emph{arc}
is an interval
component of a curve, and a \emph{loop} is a circle component.
A \emph{curve of type $+$} is a curve with only loop
components and a \emph{curve of type $-$} is a
curve with no loops parallel to the boundary.
The \emph{geometric intersection number} $i([C_1], [C_2])$ of two
curves is the
minimal number of intersections (without signs) between
representatives of the homotopy classes $[C_1]$ and $[C_2]$.

A \emph{weighted curve} is a positive linear combination of curves, where
two parallel components may be merged and their weights added. The
space of weighted simple curves on~$\Sigma$ is denoted
$\Curves(\Sigma)$. If $\Sigma$ has boundary, then we distinguish two
subsets:
\begin{itemize}
\item $\Curves^+(\Sigma)$ is the space of weighted simple curves of
  type~$+$; and
\item $\Curves^-(\Sigma)$ is the space of weighted simple curves of
  type~$-$.
\end{itemize}
\end{definition}

\begin{remark}
  As curves need not be connected, they are what other authors would
  call a \emph{multi-curve}.
\end{remark}

There are several different geometric structures one can put on a
surface.
First, we can consider a \emph{conformal} or \emph{complex}
structure~$\omega$ on~$\Sigma$, considered up to isotopy.

The next two structures deal with \emph{measured foliations} (or
equivalently
measured laminations). We always consider measured foliations up
to homotopy and Whitehead equivalence. Given a measured foliation $F$
and a
curve~$C$, we can compute $i([C],F)$, the minimal (transverse) length
of any curve isotopic to~$C$ with respect to~$F$. This is unchanged
under Whitehead equivalence, and the converse is true: two measured
foliations are Whitehead equivalent iff the transverse lengths of all
curves within an appropriate dual class (as specified below) are the same.

For a
surface~$\Sigma$ with no boundary, there is only one type of measured
foliation, and they form a finite-dimensional space
$\MF(\Sigma)$. The space of curves $\Curves(\Sigma)$ is dense in
$\MF(\Sigma)$, and
measuring lengths on $\Curves(\Sigma)$ gives an embedding
$i \co \MF(\Sigma) \hookrightarrow \RR^{\Curves(\Sigma)}$.

On a surface with boundary, measured foliations come in two
natural flavors.
$\MF^+(\Sigma)$ is the space of foliations that are parallel to
  the boundary (i.e., so the transverse length of the boundary
  is~$0$). $\Curves^+(\Sigma)$ is
  dense in $\MF^+(\Sigma)$, and measuring lengths of curves in
  $\Curves^-(\Sigma)$ gives an embedding
  $i \co \MF^+(\Sigma) \hookrightarrow \RR^{\Curves^-(\Sigma)}$.

Dually, $\MF^-(\Sigma)$ is the space of measured foliations without
boundary
  annuli. It is the closure of $\Curves^-(\Sigma)$, and there is an
  embedding $i \co \MF^-(\Sigma) \hookrightarrow\RR^{\Curves^+(\Sigma)}$.

For a closed surface~$\Sigma$, we define $\MF^+(\Sigma)$ and
$\MF^-(\Sigma)$ to be equal to $\MF(\Sigma)$.

\begin{warning}
  The set of \emph{connected} curves is not dense in $\MF^\pm(\Sigma)$ if
  $\Sigma$ has non-empty boundary.
\end{warning}

Finally, we consider quadratic differentials, in the following
variants. A \emph{quadratic differential} on~$\Sigma$ is locally of the form
  $q(z)\,(dz)^2$ with $q$ holomorphic, and determines a half-turn
  surface structure on~$\Sigma$, away from a finite number of singular
  cone points. (A \emph{half-turn surface} is a surface with
  a chart where the overlap maps are translations or rotations
  by~$\pi$.) $\Quad(\Sigma, \omega)$ is the space of quadratic
  differentials with finite area.

If $\Sigma$ has boundary, then $\Quad^{\RR}(\Sigma, \omega)$ is the
  space of quadratic differentials that are real on the boundary; this
  is isomorphic to the space of quadratic differentials on $(D\Sigma,
  D\omega)$ that are invariant with respect to the involution.
We are most interested in a further subspace.
For a Riemann surface $(\Sigma, \omega)$ with boundary,
  $\Quad^+(\Sigma, \omega)$ is the subset of $\Quad^{\RR}(\Sigma, \omega)$
  that is non-negative everywhere on each boundary
  component, or equivalently where the boundary is horizontal (rather
  than vertical). For a
  surface~$\Sigma$ with no conformal structure, $\Quad^+(\Sigma)$ is
  the space of
  pairs of a conformal structure~$\omega$ on~$\Sigma$ and a quadratic
  differential in $\Quad^+(\Sigma,\omega)$, considered up to
  isotopy. Equivalently, a point in $\Quad^+(\Sigma)$ is a
  half-turn surface structure on~$\Sigma$ with horizontal
  boundary.

  From a quadratic differential~$q\in\Quad^+(\Sigma)$, we can get
  \emph{horizontal}
  and \emph{vertical}
  measured foliations
  \begin{align*}
    q^h&\in\MF^+(\Sigma),&
    q^v&\in \MF^-(\Sigma).
  \end{align*}
  At a point $x\in\Sigma$,
  the vectors $v \in T_x\Sigma$ that are tangent to one of these
  measured foliations are those for which $q(v,v) \ge 0$ for $q^h$
  or $q(v,v) \le 0$ for $q^v$. The transverse measure
  of an arc $\gamma(t)$ is given by
\begin{align*}
  q^h(\gamma) &\coloneqq
   \int \bigl\lvert \Im \sqrt{q(\gamma'(t),\gamma'(t))} \bigl\rvert\,dt&
  q^v(\gamma) &\coloneqq
   \int \bigl\lvert \Re \sqrt{q(\gamma'(t),\gamma'(t))} \bigl\rvert\,dt.
\end{align*}

A quadratic differential is more or less the combination of two of
the other three types of data.
\begin{itemize}
\item A conformal structure and a measured foliation
  $F^+\in\MF^+(\Sigma)$ uniquely determines a quadratic differential
  in $q\in\Quad^+$ with $F^+=q^h$. This is the
  \emph{Heights Theorem} \cite{HM79:QuadDiffFol,Kerckhoff80:AsympTeich,Gardiner84:MFMinNorm,MS84:Heights,Gardiner87:Teichmuller}.
\item Let $F^+, F^-$ be a pair of measured foliations in
  $\MF^+(\Sigma)$ and $\MF^-(\Sigma)$, respectively. Then
  generically there is a unique half-turn surface
  structure $q\in\Quad^+(\Sigma)$ with the given foliations as
  horizontal and vertical foliations, respectively \cite[Theorem
  3.1]{GM91:ELTeich}.
\item On the other hand, given a conformal structure and measured
  foliation $F^-\in\MF^-(\Sigma)$, there is not always a quadratic
  differential in~$\Quad^+$ with $F^-$ as its vertical measured
  foliation. You can always double the situation and consider the
  foliation $D(F^-)$ on $D(\Sigma)$, which has an associated quadratic
  differential restricting to a unique $q \in \Quad^{\RR}(\Sigma)$, but $q$ need
  not be in
  $\Quad^+(\Sigma)$; portions of $\partial\Sigma$ might be vertical rather
  than horizontal. However, if there is such a quadratic differential,
  it is unique, as the doubling argument shows.
\end{itemize}

The types of structure on surfaces are summarized on the top
row of
Table~\ref{tab:structures}.

\begin{table}
  \centerline{\begin{tabular}{rcccc}
    \toprule
      &\multicolumn{4}{c}{Data specified}\\ \cmidrule(l) {2-5}
     & $\text{Length}/\text{Width}$&
       Length&Width&Length and width\\ \midrule
    \begin{tabular}[c]{@{}r@{}}Geometric\\Surfaces\end{tabular}
       &     \begin{tabular}[c]{@{}r@{}}Conformal\\structure\end{tabular} & $\MF^+(\Sigma)$ &
       $\MF^-(\Sigma)$&$\Quad^+(\Sigma)$\\[10pt]
    \begin{tabular}[c]{@{}r@{}}Combinatorial\\Graphs\end{tabular}
       & Elastic graph & Length graph &
       Width graph & Strip graph \\ \bottomrule
  \end{tabular}}
\medskip
  \caption{The different structures on surfaces and graphs}\label{tab:structures}
\end{table}

\subsection{Convexity}
\label{sec:convexity}

One key fact, used in Section~\ref{sec:el-surfaces}, is that there is a natural
\emph{convex structure} on $\MF^\pm(\Sigma)$.
Recall that there are (several) natural coordinates for measured foliations.
\begin{itemize}
\item On $\MF^+(\Sigma)$ of a surface with boundary, pick a maximal
  collection of non-parallel disjoint simple arcs on~$\Sigma$, and
  measure the transverse lengths of each.
\item On $\MF(\Sigma)$ of a closed surface or $\MF^-(\Sigma)$ of a
  surface with boundary, take
  Dehn-Thurston coordinates with respect to some marked pair-of-pants
  decomposition of~$\Sigma$. Normalize the Dehn-Thurston twist parameter
  so that twist~$0$ corresponds to measured foliations
  that are invariant under reversing the orientation of~$\Sigma$.
  (See, e.g., \cite{Thurston:GeomIntersect}.)
\end{itemize}
Call any of these coordinate systems \emph{canonical coordinates}.

\begin{definition}\label{def:strongly-convex}
  A function on $\MF^\pm(\Sigma)$ is
  \emph{strongly convex} if it is convex as a function on $\RR^k$ for
  each choice of canonical coordinates.
\end{definition}

This definition appears quite restrictive, since there are infinitely
many different canonical coordinates. However, such functions do
exist.
\begin{theorem}\label{thm:subadditive-convex}
  For $b \in \{+,-\}$, let $f$ be a
  function on weighted curves of type~$b$ so that
  \begin{itemize}
  \item $f$ does not increase under smoothing of essential
  crossings:
  \begin{equation}
    f\Bigl[\mfigb{curves-0}\Bigr] \ge
      f\Bigl[\mfigb{curves-1}\Bigr],\label{eq:convex-smooth}
  \end{equation}
  where the crossing is essential and the strands have the same weight, and
  \item $f$  is convex under union: for $C_1$ and $C_2$ two curves,
  \begin{equation}
    f\bigl[(C_1 \cup C_2)/2\bigr] \le \bigl(f[C_1] + f[C_2]\bigr)/2.
  \label{eq:convex-avg}
  \end{equation}
  \end{itemize}
  Then $f$ extends
  uniquely to a continuous, strongly convex
  function on $\MF^b(\Sigma)$.
\end{theorem}
An \emph{essential crossing} of a curve is (somewhat loosely) one that
cannot be removed by homotopy.
Note that $C_1 \cup C_2$ need not be simple, even if $C_1$ and~$C_2$ are.

\begin{example}
  Take any geodesic metric~$g$ on~$\Sigma$, and let $\ell_g[C]$ be the
  minimal length of a curve in the homotopy class~$[C]$ with respect
  to the metric~$g$. Then
  Equation~\eqref{eq:convex-smooth} is true, as smoothing crossings in
  the geodesic representative can only decrease length, and
  Equation~\eqref{eq:convex-avg} is true by definition (with
  equality). Thus all of these functions extend to continuous
  functions on $\MF^b(\Sigma)$. Special cases of interest include when
  \begin{itemize}
  \item the metric is a hyperbolic metric on~$\Sigma$, or
  \item the metric degenerates so that the lengths approach
    the transverse measure with respect to a measured foliation.
  \end{itemize}
  As an example of this degeneration, the function
  $F \mapsto i([C],F)$ for a fixed background curve~$C$ is a strongly
  convex function on $\MF^b(\Sigma)$.
\end{example}

Lemmas~\ref{lem:el-smoothing} and~\ref{lem:el-union} below say that
extremal length
also satisfies the same conditions.


\begin{proof}[Proof sketch of Theorem~\ref{thm:subadditive-convex}]
  Fix a canonical set of coordinates on $\MF^b(\Sigma)$.
  Given two rational measured foliations $F_0$, $F_1\in\MF^b$, let
  $F_{1/2} = (F_0 \oplus F_1)/2$ be
  the midpoint of the straight-line path between them, where $\oplus$
  means adding the chosen canonical coordinates. (This
  depends on the coordinate system.) Since they are rational, $F_0$
  and~$F_1$ are represented by weighted simple curves, as
  is~$F_{1/2}$. Analysis of the coordinates shows that
  $F_{1/2}$ is obtained from $(F_0 \cup F_1)/2$ by smoothing
  crossings (for some choice of resolution of the crossings). Thus, by
  Equations \eqref{eq:convex-smooth} and~\eqref{eq:convex-avg},
  \[
  f[F_{1/2}] \le f\bigl[(F_0\cup F_1)/2\bigr] \le \bigl(f[F_0] + f[F_1]\bigr)/2,
  \]
  which implies that $f$, when defined, is a convex
  function. But $f$ is defined on the dyadic
  rational points in $\MF^b(\Sigma)$, and
  continuity of $f$ follows from convexity \cite[Theorem
  10.1]{Rockafellar70:ConvexAnalysis}.
\end{proof}

\subsection{Graphs}
\label{sec:graphs}

\begin{figure}
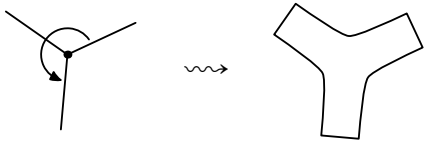

  \[
  \mfigb{surface-0} \quad\rightsquigarrow\quad \mfigb{surface-1}
  \]
  \caption{Thickening a graph into a surface, $\Gamma \rightsquigarrow
  N\Gamma$}
  \label{fig:thicken-graph}
\end{figure}

In this paper, a \emph{graph}~$\Gamma$ is a connected 1-complex, with possibly
multiple edges and self-loops. A \emph{ribbon graph} has in
addition a cyclic ordering on the ends incident to each vertex; this
gives a canonical thickening of~$\Gamma$ into an oriented surface with
boundary $N\Gamma$, as in Figure~\ref{fig:thicken-graph}.

A \emph{curve} on a graph~$\Gamma$ is a map from a $1$-manifold~$C$
into~$\Gamma$, considered up to homotopy. We do not admit
arcs here, and up to homotopy we can assume
that curves are \emph{taut}: they do not backtrack on themselves.

We can put geometric structures on graphs corresponding loosely to the
four geometric structures on surfaces above, as outlined in
Table~\ref{tab:structures}.

First, corresponding to a conformal structure, we can consider an
\emph{elastic graph}~$\Gamma$, with an
\emph{(elastic) weight} on the edges: for each edge~$e$ of~$\Gamma$,
give a positive measure which is absolutely continuous with respect to
the Lebesgue measure. Up to equivalence, we effectively just give a total
measure~$\alpha(e)$ on each edge.

There are at least two ways to interpret this
measure~$\alpha$. On one hand, we can create a \emph{rubber-band network}. For
each edge~$e$,
  take an idealized rubber band with spring constant $1/\alpha(e)$, so
  the Hooke's Law energy when the edge is stretched to length $\ell(e)$ is
  \begin{equation}\label{eq:harmonic}
  \textrm{Energy}(e) = \ell(e)^2/\alpha(e).
  \end{equation}
  Attach these rubber bands at the
  vertices. Note that, as
  for real rubber bands, a longer section of rubber band (larger
  $\alpha(e)$) is easier to stretch. Unlike for real rubber bands, the
  resting length is~$0$.

\begin{figure}
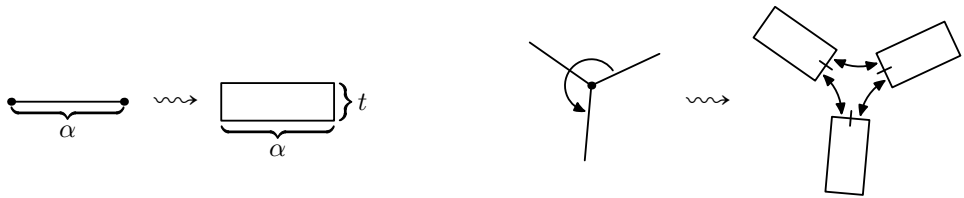

    \[
    \mfigb{surface-10}\;\rightsquigarrow\; \mfigb{surface-11}
    \qquad\qquad\qquad
    \mfigb{surface-0}\;\rightsquigarrow\; \mfigb{surface-2}
    \]
    \caption{Geometrically thickening a graph. Left: An edge of
      length~$\alpha$ is thickened to an $\alpha \times \epsilon$
      rectangle. Right: At a vertex, glue each half of the end of each
      rectangle to one of the neighbors. This is a more
      geometrically precise
      version of Figure~\ref{fig:thicken-graph}.}
    \label{fig:geom-thicken}
  \end{figure}

  Alternately, we can think of $\alpha(e)$ as defining a family of
  \emph{rectangle surfaces}: given elastic weights
  on a ribbon graph~$\Gamma$ and a constant $\epsilon > 0$, define a
  conformal surface
  $N_\epsilon \Gamma$ by taking a rectangle of size $\alpha(e) \times
  \epsilon$ for each edge~$e$ and gluing them according to the ribbon
  structure at the vertices, as indicated in
  Figure~\ref{fig:geom-thicken}.
  $\Gamma$ should be thought of as the limit of
  the Riemann surfaces $N_\epsilon\Gamma$ as $\epsilon \to 0$.

  The next geometric structure is a \emph{width graph}, which is a
  graph~$C$ with a \emph{width} on
  the edges: for each edge~$e$, there is a width $w(e)$ satisfying a
  triangle inequality at each vertex~$v$, as follows. Let
  $e_1,\dots,e_n$ be the edges incident to~$v$ (with an edge appearing
  twice if it is a self-loop). Then we require, for each $1 \le j \le
  n$,
  \begin{equation}
    \label{eq:triang}
    w(e_j) \le \sum_{\substack{1 \le i \le n\\i \ne j}} w(e_i).
  \end{equation}
  $\Wid(\Gamma)$ is the space of possible widths on an abstract
  graph~$\Gamma$.

  For a ribbon graph~$\Gamma$, there is a natural surjective map $W
  \co \MF^+(N\Gamma) \to \Wid(\Gamma)$. For $e$ an edge of~$\Gamma$,
  let $e*$ be the arc on $N\Gamma$ dual to~$e$. Then $W(F^+)$ takes
  the value $F^+(e^*)$ on the edge~$e$. When $\Gamma$ is
  trivalent (i.e., the vertices of
  $\Gamma$ have degree $\le 3$), $W$ is an isomorphism. In
  general there is ambiguity about how to glue at the vertices.

Dually, a \emph{length graph}~$K$ is a graph together with an
  assignment~$\ell$ of a non-negative length to each edge of~$\Gamma$,
  which we think of as defining a pseudo-metric on~$\Gamma$. (It is a
  pseudo-metric because points can be distance zero apart.)
  $\Len(\Gamma)$ is the space of lengths on~$\Gamma$.

  If $\Gamma$ is a ribbon graph with no valence~$1$ vertices, there is
  a natural injective
  map $L^* \co \Len(\Gamma) \to \MF^-(N\Gamma)$, defined by
  \[
  L^*(\ell) \coloneqq \sum_{e \in \Edges(\Gamma)} \ell(e) \cdot e^*
  \]
  for $\ell$ a loop on~$\Gamma$.
  This map is not surjective, but if $\Gamma$ is trivalent, $L^*$
  maps onto an open subset of $\MF^-(N\Gamma)$.

Finally, a \emph{strip graph}~$S$ is a graph in which each edge~$e$ is equipped
  with
  \begin{itemize}
  \item a length $\ell(e)$,
  \item a width $w(e)$, and
  \item an elastic weight (aspect ratio) $\alpha(e)$,
  \end{itemize}
  so that $\ell(e) = \alpha(e) w(e)$ and the widths satisfy
  the triangle inequalities, Equation~(\ref{eq:triang}).
  A strip graph has a natural underlying elastic graph $\Gamma(S)$,
  width graph $C(S)$, and length graph $K(S)$. It also has a total
  area
  \begin{equation}\label{eq:area}
  \Area(S) = \sum_{e\in \Edges(S)} \ell(e)w(e)
    = \sum_{e \in \Edges(S)} \frac{\ell(e)^2}{\alpha(e)}
    = \sum_{e \in \Edges(S)} \alpha(e)w(e)^2.
  \end{equation}
  (Compare this to elastic energy, Equation~\eqref{eq:harmonic}
  above, and extremal length on graphs, Equation~\eqref{eq:el-graph-2}
  below.) 

\begin{warning}\label{warn:elastic-vs-length}
  Be careful to distinguish between lengths and elastic weights. They
  are determined by the same data (a measure or equivalently a length
  on each edge), but they are interpreted differently. If the elastic
  weights are interpreted as defining a rubber band network, then the
  lengths can be interpreted as lengths of a system of pipes through
  which the rubber bands can be stretched. Alternatively, when the
  elastic weights are aspect ratios of rectangles (length$/$width), the
  lengths are just the length of the rectangles.
\end{warning}

As for surfaces, two of the three other types of data determine a
strip graph, except that a choice of elastic weights and lengths may
not correspond to a strip graph, as the triangle inequalities on the
widths may be violated.

\begin{remark}
  In the thickening of an elastic graph~$\Gamma$ into a surface
  $N_\epsilon \Gamma$ in Figure~\ref{fig:geom-thicken},
  the precise details
  of how you glue the rectangles at the vertices are
  irrelevant in the limit, in the following sense. If we pick two different
  ways of gluing at a vertex (e.g., gluing different proportions to
  the left and right) and get two families of surfaces
  $N_\epsilon^1\Gamma$ and $N_\epsilon^2\Gamma$, then in the limit as
  $\epsilon \to 0$ the minimal quasi-conformal constant of maps
  between $N_\epsilon^1\Gamma$ and $N_\epsilon^2\Gamma$ goes to~$1$.

  In particular, a strip graph~$S$ has a thickening $N_\epsilon S$,
  where each edge~$e$ is replaced by a rectangle of size $\ell(e)
  \times \epsilon w(e)$, and the rectangles are glued in a
  width-preserving way at the vertices. (This gluing is possible by
  the triangle inequalities.) If $\Gamma(S)$ is the underlying elastic
  graph of~$S$, then for $\epsilon \ll 1$ the surfaces
  $N_\epsilon S$ and $N_\epsilon \Gamma(S)$ are nearly conformally
  equivalent. See Proposition~\ref{prop:el-graphs-surfaces} for
  some related estimates.
\end{remark}

\subsection{Extensions}
\label{sec:extensions}

The graphs $\Gamma$ we are considering in this paper are usually
topological spines of the corresponding surface $N\Gamma$ (i.e.,
$N\Gamma$ deformation retracts onto~$\Gamma$). More generally, we may
consider a graph~$\Gamma$ and surface~$\Sigma$ with a
$\pi_1$-surjective embedding $\Gamma \hookrightarrow \Sigma$, or more
generally still a graph~$\Gamma$, group~$G$, and surjective map
$\pi_1(\Gamma) \to G$ (i.e., a generating graph for~$G$). Much of the
theory extends to this case.

It is also sometimes convenient to generalize from surfaces to
orbifolds. In the setting of groups, this means considering
maps from $\pi_1(\Gamma)$ to the orbifold fundamental group
$\pi_1^{\textrm{orb}}(\Sigma)$.
In the setting of graphs embedded in surfaces, we consider graphs with
marked
points that are required to be mapped to marked (orbifold) points. In
particular,
this is required for a proper statement of
Theorem~\ref{thm:detect-rational} if we drop the restriction that
there be a branch point in each cycle of~$P$.


\section{Harmonic maps and Dirichlet energy}
\label{sec:harmonic}
We now turn to harmonic maps, maps that minimize some form of
Dirichlet energy.
\subsection{Harmonic maps from surfaces}
\label{sec:harmonic-surfaces}
Given a conformal Riemann surface $(\Sigma, \omega)$, a length
graph~$(K, \ell)$, and a Lipschitz map $f \co \Sigma \to K$, the
\emph{Dirichlet energy} of $f$ is
\begin{equation}
  \label{eq:dirichlet-surface}
  \Dir(f) \coloneqq \int_\Sigma \abs{\nabla f(z)}^2\,dA.
\end{equation}
Here, we have picked an (arbitrary) Riemannian metric~$g$ in the given
conformal class, $dA$ is the area measure with respect to~$g$, and
$\abs{\nabla f(z)}$ is defined to be the best Lipschitz constant of $f$
at~$z$ with respect to~$g$ and~$\ell$. (This agrees with the usual norm of
the gradient when $f(z)$ is in the interior of an edge
of~$K$ and $f$~is differentiable.)

We also define $\Dir[f]$ to be the lowest energy in the homotopy class of~$f$:
\[
\Dir[f] \coloneq \inf_{g \in [f]} \Dir(g).
\]
This optimum is achieved \cite[Theorem 11.1]{EF01:HarmonicPolyhedra}.
The optimizing functions are called \emph{harmonic maps} and define
canonical quadratic differentials in $\Quad^+(\Sigma)$.

\begin{warning}\label{warn:harmonic-vs-quad-diff}
  The connection between harmonic maps and quadratic differentials is
  not as tight as you might expect. Given a closed surface~$\Sigma$
  and a weighted simple curve~$C$, we can construct a length
  graph~$K(C)$, with vertices the connected components of $\Sigma
  \setminus C$ and edges given by components of~$C$, with lengths of
  edges of $K(C)$ given by weights on~$C$. There is a natural homotopy
  class of maps $\Sigma \to K(C)$. The Dirichlet minimizer in
  this homotopy class sometimes, but not always, recovers the
  Jenkins-Strebel quadratic differential with vertical foliation given
  by~$C$.

  A similar construction using maps from the universal cover
  of~$\Sigma$ to an $\RR$-tree
  does recover the quadratic differential with given vertical
  foliation. Wolf used this to
  give an alternate proof of the Heights Theorem
  \cite{Wolf98:MeasuredFoliations}.

  Alternatively, we could generalize the target space, and allow $K$
  to be a general non-positively curved polyhedral complex. In this
  context, we can take a complex $K_\epsilon(C)$ made of very thin
  tubes around the edges around the graph $K(C)$ defined above. Then
  there is a
  closer connection between harmonic maps to $K_\epsilon(C)$ and the
  Jenkins-Strebel quadratic differential. We do not pursue this
  further here.
\end{warning}

\subsection{Harmonic maps from graphs}
\label{sec:harmonic-graphs}

Given an elastic graph $(\Gamma, \alpha)$ and a length graph $(K,
\ell)$, the \emph{Dirichlet energy} of a Lipschitz map $f \co \Gamma \to K$ is
\begin{align}
  \label{eq:dirichlet-graph}
  \Dir(f) &\coloneqq \int_\Gamma \abs{f'(x)}^2\,dx = \int_K \Fill_f(y)\,dy\\
  \Dir[f] &\coloneqq \inf_{g \in [f]} \Dir(g).
\end{align}
Here, $\abs{f'(x)}$ is the derivative of $f$ with respect to the
natural coordinates given by $\alpha$ on~$\Gamma$ and $\ell$
on~$K$. This derivative is not defined at vertices, but these points
are negligible. Alternatively, $\abs{f'(x)}$ is the best Lipschitz
constant of $f$ at~$x$. $\Fill_f(y)$ is the \emph{filling function} of
$f$ at $y\in K$, the sum of derivatives at preimages:
\begin{equation}\label{eq:fill-def}
\Fill_f(y) \coloneqq \sum_{f(x) = y} \abs{f'(x)}.
\end{equation}
The two expressions for $\Dir(f)$ are related by an easy change of
variables.

Minimizers of Dirichlet energy are harmonic
functions in the following sense.

\begin{definition}\label{def:harmonic-graphs}
  A function $f \co \Gamma \to K$ from an elastic graph~$\Gamma$ to a
  length graph~$K$ is \emph{harmonic} if the following conditions are
  satisfied.
  \begin{enumerate}
  \item The map $f$ is piecewise linear.
  \item The map $f$ does not backtrack (i.e., $f$ is locally injective on
    each edge of $\Gamma$).
  \item The derivative $\abs{f'(x)}$ is constant on the edges of~$\Gamma$ when
    defined.\label{item:deriv-const}
    As a result, for $e$ an edge
    of~$\Gamma$ we may write $\abs{f'(e)}$ for the common value at any
    point on the edge.
  \item\label{item:cond-edge-balance} If a vertex~$v$ of~$\Gamma$ maps
    to the interior of an edge
    of~$K$, with edges $e_1,\dots,e_k$ incident on the left and edges
    $e_{k+1},\dots,e_{k+l}$ incident on the right, then
    \begin{equation}
      \label{eq:vertex-edge-balance}
      \sum_{i=1}^k \abs{f'(e_i)} = \sum_{i=1}^l \abs{f'(e_{i+k})}.
    \end{equation}
    In particular, for each edge $E$ of~$K$, the filling function
    $\Fill_f(y)$ is constant on~$E$, and we may write
    $\Fill_f(E)$.
  \item\label{item:cond-vertex-balance} If a vertex~$v$ of~$\Gamma$
    maps to a vertex~$w$ of~$K$, let
    $E_1,\dots,E_k$ be the germs of edges of~$K$ incident to~$w$.
    Then we have a vertex balancing condition: for $i =
    1,\dots,k$,
    \begin{equation}
      \sum_{\substack{\text{$e$ incident to $v$}\\e \mapsto E_i}} \abs{f'(e)} \le
        \sum_{\substack{1 \le j \le k\\j \ne i}} \,\sum_{e \mapsto E_j} \abs{f'(e)}.
      \label{eq:vertex-vertex-balance}
    \end{equation}
    Here, the notation ``$\sum_{e \mapsto E_i}$'' means the sum over
    all germs of edges of~$\Gamma$ that locally map to~$E_i$.
  \end{enumerate}
\end{definition}

\begin{citethm}\label{thm:harmonic-graphs}
  A function $f \co \Gamma \to K$ is a local minimum for the Dirichlet
  energy within its homotopy class iff it is harmonic. Every local
  minimum is also a global minimum.
\end{citethm}

From the rubber bands point of view, the intuition behind
Definition~\ref{def:harmonic-graphs} is that $2\abs{f'(x)}$ is the
\emph{tension} in the rubber band, which is constant along the
edge. The net force on a vertex mapping
to the interior of an edge must be zero; this is
condition~(\ref{item:cond-edge-balance}). For a vertex~$v$ mapping to
a vertex~$w$, the net force pulling $v$ in any one direction away
from~$w$ cannot be too large; this is
condition~\eqref{item:cond-vertex-balance}.
Condition~\eqref{item:cond-edge-balance} can be thought of as the
special case of condition~\eqref{item:cond-vertex-balance} when $w$ has
only two incident edges.

\begin{figure}
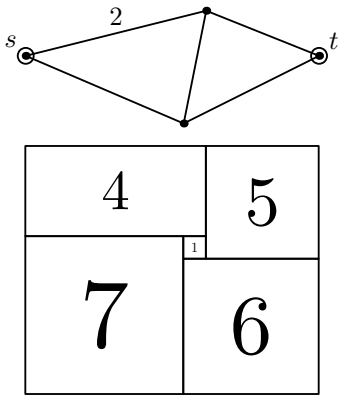

  \[
  \mfig{squares-1}
  \]
  \[
  \mfig{squares-0}
  \]
  \caption{A simple marked elastic graph and the associated tiling by
    rectangles. All elastic weights on the graph are~$1$ except for
    one edge, which has weight~$2$. On the rectangle tiling, we have shown the
    widths of the rectangles, which are $\abs{f'(e)}$ in the
    associated harmonic map to the interval.}
  \label{fig:rect-tiling}
\end{figure}

Conditions~\eqref{item:cond-edge-balance}
and~\eqref{item:cond-vertex-balance} also imply that the derivatives
$\abs{f'(e)}$ form a valid width structure on~$\Gamma$ and that the
filling functions $\Fill_f(E)$ form a width structure on~$K$.
From the rectangular surface point of view,
a harmonic map gives a tiling of~$K$ with
rectangles, with each edge~$e$ of~$\Gamma$ contributing a rectangle of
aspect ratio $\alpha(e)$. To illustrate this, consider the case of
marked planar graphs. Let $\Gamma$ be a planar graph with two
distinguished vertices $s$ and~$t$ bordering the infinite face, and
consider the minimizer of the Dirichlet energy among maps
\[
(\Gamma, s, t) \longrightarrow ([0,1], \{0\}, \{1\})
\]
mapping $\Gamma$ to the interval, taking $s$ to~$0$ and $t$
to~$1$. Brooks, Smith, Stone, and Tutte showed that this minimizer
gives a rectangle packing of a rectangle
\cite{BSST40:DissectionRects}. See Figure~\ref{fig:rect-tiling} for an
example.

\begin{remark}\label{rem:not-electrical}
  In the context of rectangle packings, it is traditional to use the
  language of
  resistor networks rather than elastic or rubber-band networks. The
  Dirichlet energy in these two settings is the same for maps to an
  interval, but for
  general target graphs rubber bands are more flexible, because of the
  lack of orientations. Where rubber band optimizers are related to
  homotopy of maps from~$\Gamma$ to a target graph and
  holomorphic quadratic differentials, electrical current or resistor
  network optimizers are related to homology of~$\Gamma$ and to
  holomorphic differentials (not quadratic differentials).
\end{remark}

\subsection{Relation to Lipschitz energy}
\label{sec:lipschitz-energy}

Another natural norm for maps $\phi \co K_1 \to K_2$
between metric graphs is the \emph{Lipschitz energy}
\begin{align}
  \label{eq:lipschitz}
  \Lip(\phi) &\coloneqq \esssup_{x \in K_1} \abs{\phi'(x)}\\
  \Lip[\phi] &\coloneqq \inf_{\psi \in [\phi]} \Lip(\psi).
\end{align}
Compared to Equation~\eqref{eq:dirichlet-graph}, we take the
$L^\infty$ norm of $\abs{\phi'(x)}$ rather than the $L^2$ norm.

Let us note that
\begin{itemize}
\item $\Lip(\phi) < 1$ iff $\phi$ is distance-decreasing, and
\item the Lipschitz energy is
  \emph{dynamical}, in the sense that it behaves well under
  composition:
  \begin{equation*}
    \Lip[\phi \circ \psi] \le \Lip[\phi] \Lip[\psi].
  \end{equation*}
\end{itemize}
The Lipschitz energy should be thought of as an invariant of a
map between two length graphs, not as an invariant of a map from an
elastic graph to a length graph.

It follows immediately from the definitions that $\Lip$ is
sub-multiplicative with respect to both length and Dirichlet energy:
for $C$ a curve on $K_1$ and $f \co \Gamma \to K_1$ a map from an
elastic graph,
\begin{align*}
  \ell_{K_2}[\phi_* C] &\le \Lip[\phi]\ell_{K_1}[C]\\
  \Dir[\phi\circ f] &\le \Lip[\phi]^2\Dir[f].
\end{align*}
These are both special cases of
Equation~\eqref{eq:norm-submultiplicative}, and in both cases the
inequalities are tight. For instance, we have \cite[Proposition
3.11]{FM11:MetricOutSpace}
\begin{equation}\label{eq:lip-ratio}
  \Lip[\phi] = \sup_{\text{$C$ curve on $K_1$}}
    \frac{\ell_{K_2}[\phi_* C]}{\ell_{K_1}[C]}.
\end{equation}

\subsection{Computing Dirichlet energy}
\label{sec:compute-Dirichlet}
Given a homotopy class $[f] \co \Gamma \to K$ of maps between
an elastic graph~$\Gamma$ and a metric
graph~$K$, how can one find a representative $g \in [f]$ that
minimizes Dirichlet energy? By Theorem~\ref{thm:harmonic-graphs}, this
is the same as finding a
harmonic map in the homotopy class.

As with harmonic maps in other settings,
this is easy to compute, with at least two reasonable approaches:
\begin{enumerate}
\item \emph{Repeated averaging}, iteratively moving the image of each
  vertex of~$\Gamma$ to the weighted average of its neighbors in the
  universal cover $\wt K$ of~$K$. The \emph{average}
  of a non-empty set~$S$ of points on a metric tree is the (unique) point~$x$
  that minimizes the sum of squares of distances from $x$ to points
  in~$S$. Note that this average can generically be on a vertex.
\item \emph{Linear} or \emph{convex quadratic programming}, based on
  the observation that once the combinatorics of the map are fixed,
  specifying which vertices of~$\Gamma$ go to which vertices or edges
  of~$K$, the Dirichlet energy is a convex, quadratic function of the
  positions along the edges, and thus the minimum can be found by
  solving linear equations. Here, one starts by guessing some
  combinatorics, and updating the combinatorics if it turns out not to
  be optimal.
\end{enumerate}
The second approach should usually be faster, but also requires some
more care in changing the combinatorics. The harmonic representative
of $[f]$
is not unique in general; however, the set of harmonic representatives
forms a convex set in a suitable sense.


\section{Extremal length}
\label{sec:extremal-length}

\subsection{Extremal length on surfaces}
\label{sec:el-surfaces}

Given a Riemann surface $(\Sigma, \omega)$ and a curve~$C$ on~$\Sigma$,
recall that the \emph{extremal length} of~$C$ on~$\Sigma$ is (omitting
analytic details)
\begin{equation}
  \label{eq:el-surface}
  \EL_\omega[C] \coloneqq
    \sup_{\rho: \Sigma \to \RR_{\ge 0}} \frac{\ell_{\rho
        g}[C]^2}{\Area_{\rho g}(\Sigma)},
\end{equation}
where we make the following definitions.
\begin{itemize}
\item The metric $g$ is an arbitrary metric in the conformal class~$\omega$.
\item The metric $\rho g$ is the metric~$g$ scaled by the conformal
  factor~$\rho$. (This may be a pseudo-metric.)
\item The number $\ell_{\rho g}[C]$ is the minimal length of any
  element of the homotopy class~$[C]$ in the
  metric~$\rho g$.
\item The number $\Area_{\rho g}(\Sigma)$ is the total area of~$\Sigma$ with respect
  to~$\rho g$.
\end{itemize}
Observe that scaling~$\rho$ by a
global constant does not change the supremand.

We can extend Equation~\eqref{eq:el-surface} to allow for weighted
curves: for $C = \sum_i a_i C_i$ a weighted curve, define $\ell[C]
\coloneqq \sum_i a_i \ell[C_i]$, and use
Equation~\eqref{eq:el-surface} as before. (Note that this agrees with
the earlier definition for integral weights.)

The following results are standard.

\begin{citethm}[Jenkins-Strebel]
  If $C$ is a weighted simple curve, then the supremum in the
  definition of $\EL[C]$ is achieved, and the optimal metric is the
  metric on a
  half-turn surface associated to a quadratic differential with
  horizontal foliation equal to~$C$.
\end{citethm}

\begin{lemma}\label{lem:el-cover}
  Let $\pi \co \wt \Sigma \to \Sigma$ be a covering map of degree~$d$. For $C$
  a weighted curve on~$\Sigma$, define $\pi^{-1} C$ to be the
  inverse image of~$C$, with the same weights. Then $\EL[\pi^{-1} C] = d\EL[C]$.
\end{lemma}

\begin{lemma}\label{lem:el-quadratic}
  If $a \in \RR_{> 0}$ is a global weighting factor, then
  \[
  \EL[a C] = a^2 \EL[C].
  \]
\end{lemma}

We now turn to properties related to convexity, as in
Section~\ref{sec:convexity}.
The following two properties follow from elementary arguments.

\begin{lemma}\label{lem:el-smoothing}
  $\EL$ does not increase under smoothing of essential crossings:
  \[
  \EL\Bigl[\mfigb{curves-0}\Bigr] \ge \EL\Bigl[\mfigb{curves-1}\Bigr].
  \]
  Here, the two sides show a local picture of unweighted curves, or
  weighted curves with equal weights on the two local strands.
\end{lemma}

\begin{lemma}\label{lem:el-union}
  If $C_1$ and $C_2$ are two weighted multi-curves, then
  \[
  \EL[C_1 \cup C_2] \le 2(\EL[C_1] + \EL[C_2]),
  \]
  where we keep the weights on each component of $C_1$ and~$C_2$.
  More generally, for $0 \le t \le 1$,
  \[
  \EL[t C_1 \cup (1-t) C_2] \le t \EL[C_1] + (1-t)\EL[C_2].
  \]
\end{lemma}

\begin{remark}
  In Lemma~\ref{lem:el-union}, we are taking the union of curves, not
  the union of path families, as sometimes appears in the theory of
  extremal length.
\end{remark}

\begin{corollary}\label{cor:el-continuous}
  $\EL$ extends uniquely to a continuous, strongly convex function on
  $\MF^\pm(\Sigma)$.
\end{corollary}

\begin{proof}
  Follows from Theorem~\ref{thm:subadditive-convex}.
\end{proof}

\begin{remark}
  Presumably the techniques in the proof of
  Corollary~\ref{cor:el-continuous} can be extended to prove the rest
  of the Heights Theorem, in particular that the associated quadratic
  differential varies continuously as a function of the measured
  foliation.
\end{remark}

\subsection{Extremal length on graphs}
\label{sec:el-graphs}

By analogy with Equation~\eqref{eq:el-surface}, for $C$ a weighted
curve on an elastic graph $(\Gamma, \alpha)$,
define
\begin{equation}
  \label{eq:el-graph}
  \EL_\alpha[C] \coloneqq
    \sup_{\rho: \Edges(\Gamma) \to \RR_{\ge 0}}
      \frac{\ell_{\rho \alpha}[C]^2}{\Area_{\rho\alpha}(\Gamma)},
\end{equation}
where
\begin{itemize}
\item The length metric $\rho\alpha$ on~$\Gamma$ gives edge~$e$ the length $\rho(e)\alpha(e)$.
\item The number $\ell_{\rho\alpha}[C]$ is the length of $C$ with respect to
  $\rho\alpha$, i.e.,
  \begin{equation}
    \label{eq:len-graph}
    \ell_{\rho\alpha}[C] \coloneqq \sum_{e \in \Edges(\Gamma)} n_C(e) \rho(e)\alpha(e),
  \end{equation}
  where $n_C(e)$ is the weighted number of times that $C$ runs
  over~$e$.
\item $\Area_{\rho\alpha}(\Gamma)$ is the ``area'' of $\Gamma$ with
  respect to $\rho\alpha$, defined to be
  \begin{equation}
    \label{eq:area-graph}
    \Area_{\rho\alpha}(\Gamma) \coloneqq \sum_{e \in \Edges(\Gamma)} \rho(e)^2 \alpha(e).
  \end{equation}
  (The intuition is that each edge is turned into a rectangle of width
  proportional to $\rho(e)$ and aspect ratio $\alpha(e)$, and thus
  area $\rho(e)\alpha(e)$.)
\end{itemize}

In fact the supremum in Equation~\eqref{eq:el-graph} is easy to
do. The optimum has $\rho(e)$ proportional to $n_C(e)$ and so
\begin{equation}
  \label{eq:el-graph-2}
  \EL_\alpha[C] = \sum_{e \in \Edges(\Gamma)} n_C(e)^2 \alpha(e).
\end{equation}
This formula extends immediately to a function on $\Wid(\Gamma)$, and
satisfies Lemmas~\ref{lem:el-cover}, \ref{lem:el-quadratic},
and~\ref{lem:el-union}. For
Lemma~\ref{lem:el-smoothing}, we need to pick a ribbon structure
on~$\Gamma$ in order to make sense of ``essential'' crossings. With
any such choice, Lemma~\ref{lem:el-smoothing} is true.

\begin{corollary}
  Let $(\Gamma,\alpha)$ be an elastic spine for a
  surface~$\Sigma$. Then the function $C \mapsto \EL_\alpha[C]$ on
  $\Curves^+(\Sigma)$ extends uniquely to a strongly convex function
  on $\MF^+(\Sigma)$.
\end{corollary}

\subsection{Relating graphs and surfaces}
\label{sec:el-graphs-surfaces}

We can now give a concrete relation between an elastic graph
$(\Gamma, \alpha)$ and the associated family of conformal surfaces
$N_\epsilon\Gamma$. Write $\EL[C;\Gamma]$ for extremal length with
respect to the elastic graph~$\Gamma$, and $\EL[C;\Sigma]$ for
extremal length with respect to the conformal surface~$\Sigma$.

\begin{proposition}\label{prop:el-graphs-surfaces}
  Let $(\Gamma,\alpha)$ be an elastic ribbon graph with trivalent vertices,
  and let $m = \min_e \alpha(e)$ be the smallest weight of any edge
  in~$\Gamma$.
  Then, for $t < m/2$ and $C$ any measured foliation on~$\Gamma$, we have
  \[
  \EL[C; \Gamma] \le t \EL[C; N_t \Gamma] \le
    \EL[C; \Gamma]\cdot (1 + 8t/m).
  \]
\end{proposition}

The proof involves finding, on the one hand, embeddings of
sufficiently thick annuli into $N_t \Gamma$, and, on the other hand,
suitable test functions~$\rho$ on $N_t \Gamma$ in
Equation~\eqref{eq:el-surface}.

\begin{remark}
  The restriction to trivalent graphs in
  Proposition~\ref{prop:el-graphs-surfaces} can be removed. Note that
  the estimate depends only on the local geometry of~$\Gamma$, and
  thus is unchanged under covers.
\end{remark}

\subsection{Duality with Dirichlet energy}
\label{sec:duality-dirichlet}

As mentioned earlier, extremal length is in some sense dual to
Dirichlet energy. More precisely, we have the following.

\begin{proposition}[Sub-multiplicative]\label{prop:dir-el-submult}
  Let $C$ be a curve in an elastic graph~$\Gamma$, and let $f \co
  \Gamma \to K$ be a harmonic map to a length graph. Then
  \[
  \ell[f_* C]^2 \le \Dir[f] \EL[C].
  \]
\end{proposition}

\begin{proof}[Proof sketch]
  This is basically the definition of $\EL$ in
  Equation~\eqref{eq:el-graph}, noting the similarities between the
  area in Equation~\eqref{eq:area-graph} and Dirichlet energy,
  Equation~\eqref{eq:dirichlet-graph}.
\end{proof}

\begin{proposition} [Duality 1]\label{prop:dir-el-dual}
  Let $f \co \Gamma \to K$ be a harmonic map from an elastic graph to
  a length graph. Then there is a sequence of weighted curves~$C_i$ in~$\Gamma$
  so that for all $\epsilon > 0$ there is an $i$ so that
  \[
  \ell[f_*C_i]^2 \ge (1-\epsilon)\Dir[f] \EL[C_i].
  \]
\end{proposition}

\begin{proof}[Proof sketch]
  Take weighted curves $C_i$ so that $n_{C_i}(e)$ approximates $\abs{f'(e)}$.
\end{proof}

\begin{proposition}[Duality 2]\label{prop:el-dir-dual}
  Let $C$ be a curve in an elastic graph~$\Gamma$. Then there is a
  length graph~$K$ and a harmonic map $f \co \Gamma \to K$ so that
  \[
  \ell[f_*C]^2 = \Dir[f] \EL[C].
  \]
\end{proposition}

\begin{proof}[Proof sketch]
  Take $K$ to be $\Gamma$ with edge lengths $\alpha(e) n_C(e)$.
\end{proof}

The situation is less satisfactory for
surfaces. Sub\hyp multiplicativity in the sense of
Proposition~\ref{prop:dir-el-submult} is true (at least when the
curve~$C$ is embedded), as is Proposition~\ref{prop:dir-el-dual}.
Issues related to Warning~\ref{warn:harmonic-vs-quad-diff} make an
analogue of Proposition~\ref{prop:el-dir-dual} more delicate, although
it is true if we allow non-positively curved polyhedral complexes as the
target space.

\begin{remark}
  The definition of extremal length on graphs in Equation
  \eqref{eq:el-graph} is somewhat backwards, in that it is a
  supremum of a ratio of energies over all metrics (or equivalently
  over all maps to length graphs). By analogy with Dirichlet energy
  (Equation \eqref{eq:dirichlet-graph}), it would be better to define
  the energy of a homotopy class as an infimum of some energy
  functional. Indeed, we could take Equation \eqref{eq:el-graph-2} as
  the primary definition.

  This remark applies to extremal length on surfaces
  (Equation \eqref{eq:el-surface}) as well: it might be better to take
  a different definition as primary. Namely, recall that a conformal
  annulus~$A$ has an extremal length $\EL(A)$, which we may define as
  the inverse of the modulus. Then
  extremal length of a weighted multi-curve $C = \sum_i a_i C_i$ can be
  alternately defined as
  \begin{equation}
    \label{eq:el-surface-2}
    \EL[C] = \inf_{A_i} \,\sum_i a_i \EL(A_i),
  \end{equation}
  where the infimum runs over all disjoint embeddings of conformal annuli
  $A_i$ with core curves homotopic to $C_i$.
\end{remark}

\subsection{Extensions}
\label{sec:extensions-el}

There are several ways in which we can extend these notions of
extremal length. First, we can consider graphs that are embedded in a
surface, not necessarily as a spine.

\begin{definition}
  A graph~$\Gamma$ embedded in a surface~$\Sigma$ is \emph{filling} if
  each component of $\Sigma \setminus \Gamma$ is a disk or an annulus
  on the boundary of~$\Sigma$.
\end{definition}

\begin{definition}\label{def:el-graph-surf}
  Let $(K, \ell)$ a length graph with a filling embedding~$\phi$ in a
  surface~$\Sigma$. For $C$ a curve in~$\Sigma$, define
  \[
  \ell_K[C] \coloneqq
    \inf_{\substack{\text{$C'$ curve on $K$}\\\phi(C') \sim C}} \ell[C'].
  \]
  Now for $(\Gamma,\alpha)$ a filling elastic graph in~$\Sigma$,
  define $\EL_{\Gamma,\alpha}[C]$ by
  \begin{equation*}
  \EL_{\Gamma,\alpha}[C] \coloneqq
    \sup_{\rho: \Edges(\Gamma) \to \RR_{\ge 0}}
      \frac{\ell_{\Gamma,\rho \alpha}[C]^2}{\Area_{\rho\alpha}(\Gamma)}.
  \end{equation*}
\end{definition}

  This is just like Equation~\eqref{eq:el-graph}, except that we
  consider homotopy classes in~$\Sigma$ rather than in~$\Gamma$.
A similar notion of extremal length was considered by Duffin
\cite{Duffin62:ELNetwork}, in the context of electrical networks and
graphs with two marked points (as in Section~\ref{sec:harmonic-graphs}).

  The optimization in Definition \ref{def:el-graph-surf} is no longer
  as easy, and the analogue of Equation~\eqref{eq:el-graph-2} is more
  awkward to state. The result of the optimization gives a
  rectangular tiling of the surface with aspect ratios given
  by~$\alpha$, analogous to Figure~\ref{fig:rect-tiling}.

More generally, we can consider graphs generating a group.

\begin{definition}\label{def:el-group}
  For $(\Gamma, \alpha)$ an elastic graph, $\phi \co \pi_1(\Gamma)
  \twoheadrightarrow G$ a surjective homomorphism onto a group~$G$, and $[g]$ a
  conjugacy class in~$G$, define
  \begin{align*}
    \ell_{\rho\alpha}[g] &\coloneqq
      \inf_{\substack{\text{$C$ curve on $\Gamma$}\\\phi[C] = [g]}}
        \ell_{\rho\alpha}[C]\\
    \EL'_\alpha[g] &\coloneqq
    \sup_{\rho: \Edges(\Gamma) \to \RR_{\ge 0}}
      \frac{\ell_{\rho \alpha}[g]^2}{\Area_{\rho\alpha}(\Gamma)}\\
    \EL_\alpha[g] &\coloneqq
      \lim_{n \to \infty} \EL'_\alpha[g^n]/n^2.
  \end{align*}
\end{definition}

There are also other models for defining a combinatorial length of
curves on a graph. Notably, Schramm \cite{Schramm93:Square} and Canon, Floyd, and
Parry \cite{CFP94:SquaringRectangles} define a model where the length of a curve in a graph
is determined by the vertices that it passes through, rather than the
edges that it crosses over (as in this paper). More generally, one can
consider \emph{shinglings} of a graph or surface, decompositions of the space
into a finite number of overlapping open sets. The combinatorial
length of a curve is then determined by which shingles it passes
through. The edge model that is the main focus of this paper comes
from taking shingles that are neighborhoods of the edges (overlapping
at the vertices), while the vertex model of
Schramm--Cannon--Floyd--Perry comes from taking shingles that are
neighborhoods of the vertices (overlapping at the centers of edges).

For any of these notions of length, one can define a notion of
extremal length using Equation~\eqref{eq:el-graph}.

More generally, one may instead attempt to characterize which extremal
length functions can appear. There are several natural sources of
an ``extremal length'' function on $\MF^+(\Sigma)$:
\begin{enumerate}
\item A conformal structure $\omega$ on $\Sigma$ gives the usual
  notion of extremal length.\label{item:el-surf}
\item If $\Sigma'$ is another surface and $\phi \co \Sigma'
  \hookrightarrow \Sigma$ is a filling embedding (an embedding
  for which the complementary regions are disks or annuli), then a
  conformal structure $\omega'$ on~$\Sigma'$ gives a notion of
  extremal length on $\Sigma$, defined analogously to
  Definition~\ref{def:el-graph-surf}.\label{item:el-surf-surf}
\item An elastic graph spine $(\Gamma,\alpha)$ for $\Sigma$ gives a
  notion of extremal length by Equations~\eqref{eq:el-graph}
  or~\eqref{eq:el-graph-2}.\label{item:el-graph}
\item A filling elastic graph in~$\Sigma$ gives a notion of
  extremal length by Definition~\ref{def:el-graph-surf}.\label{item:el-graph-surf}
\item Finally, a shingling of~$\Sigma$ gives yet another notion of
  extremal length.\label{item:el-shingle}
\end{enumerate}
All of these notions of extremal length give a function on
$\Curves^+(\Sigma)$ that
\begin{itemize}
\item is positive,
\item is homogeneous quadratic, in the sense of Lemma~\ref{lem:el-quadratic},
\item does not increase under smoothing, in the sense of
  Lemma~\ref{lem:el-smoothing},
\item is sub-additive under union, in the sense of
  Lemma~\ref{lem:el-union}, and therefore
\item extends to a strongly convex function on $\MF^+(\Sigma)$, by
  Theorem~\ref{thm:subadditive-convex}.
\end{itemize}
As a result of convexity, we can think of $\EL$ as a kind of ``norm''
on $\MF^+(\Sigma)$.
\begin{problem}\label{prob:char-el}
  Which functions $\EL\co \MF^+(\Sigma) \to \RR$ can arise from the
  constructions above?
\end{problem}
The properties above are some restrictions, but are probably not a
complete list.

Some of these notions of extremal length subsume the others: Extremal
lengths from notion \eqref{item:el-shingle} include extremal
lengths from notion~\eqref{item:el-graph-surf} by a direct
construction. By Proposition~\ref{prop:el-graphs-surfaces},
notion~\eqref{item:el-surf-surf} is dense in
notion~\eqref{item:el-graph-surf}  (up to
scale). Notion~\eqref{item:el-graph-surf}
naturally includes notion~\eqref{item:el-graph}, and by taking the
graph to be the edges of a triangulation we can see that
notion~\eqref{item:el-graph-surf} is dense in
notions~\eqref{item:el-surf} and~\eqref{item:el-surf-surf}.

\begin{problem}\label{prob:el-emb-graphs}
  The definition of extremal length in surfaces,
  Equation~\eqref{eq:el-surface}, extends to (width) graphs embedded
  in~$\Sigma$, rather than just curves. What do the resulting optimal
  metrics look like?
\end{problem}


\section{Stretch factors and embedding energy}
\label{sec:stretch-factors}

Now we study how extremal length and Dirichlet energy change under
maps between surfaces and graphs. This material is developed more
fully in a sequel paper with Pilgrim and Kahn
\cite{KPT15:EmbeddingEL}.

\subsection{Stretch factors for surfaces}
\label{sec:stretch-surface}

Let $\phi \co \Sigma_1 \to \Sigma_2$ be a topological embedding of
surfaces. Then composing with $\phi$ and deleting
null-homotopic components induces a natural map $\phi_*\co\Curves^+(\Sigma_1)
\to \Curves^+(\Sigma_2)$. (This pushforward does not work on
$\Curves^-(\Sigma_1)$.)

\begin{warning}\label{warn:map-mf-cont}
  The map $\phi_*$ does \emph{not} generally extend to a continuous map
  from $\MF^+(\Sigma_1)$ to $\MF^+(\Sigma_2)$.
\end{warning}

Now suppose $\Sigma_1$ and~$\Sigma_2$ have conformal structures
$\omega_1$ and~$\omega_2$, respectively. (The map $\phi$ need not
respect the conformal structures.) 

\begin{definition}\label{def:sf-surface}
  In the above setting, the \emph{stretch factor} of $\phi$ is
  \begin{equation}
    \label{eq:sf-surface}
    \SF[\phi] \coloneqq
      \sup_{C \in \Curves^+(\Sigma)} \frac{\EL_{\omega_2}[\phi_* C]}{\EL_{\omega_1}[C]}.
  \end{equation}
  This depends only on the homotopy class of~$\phi$.
\end{definition}

It follows from the definition that $\SF$ behaves well
under composition.

\begin{proposition}\label{prop:SF-compose}
  If $f\co \Sigma_1 \hookrightarrow \Sigma_2$ and $g \co \Sigma_2
  \hookrightarrow \Sigma_3$ are two topological embeddings of
  conformal surfaces, then
  \[
  \SF[f \circ g] \le \SF[f] \cdot\SF[g].
  \]
\end{proposition}

\begin{definition}
  A conformal embedding $\phi \co \Sigma_1
  \hookrightarrow \Sigma_2$ is \emph{strict} if, in each component
  of~$\Sigma_2$, there is a non-empty open
  subset in the complement of the image.
\end{definition}

\begin{theorem}[essentially Ioffe]\label{thm:SF-Teich}
  If $(\Sigma_1,\omega_1)$ and $(\Sigma_2,\omega_2)$ are two
  conformal surfaces and $\phi\co \Sigma_1 \to
  \Sigma_2$ is a homeomorphism so that there is no strict conformal
  embedding
  in the homotopy class~$[\phi]$, let $K$ be the lowest constant so that
  there is a $K$--quasi-conformal map in the homotopy
  class~$[\phi]$. Then
  \begin{equation*}
    \SF[\phi] = K.
  \end{equation*}
\end{theorem}

\begin{proof}
  This is very close to the main theorem of
  \cite{Ioffe75:QCImbedding}. In that paper, Ioffe proves that if
  there is no conformal embedding in $[\phi]$, there
  are canonical quadratic differentials $q_i \in \Quad^+(\Sigma_i)$
  and a $K$--quasi-conformal representative for~$\phi$ that uniformly
  stretches $q_1$
  to~$q_2$. Suitably approximating the horizontal foliations
  of $Q_1$ and~$Q_2$ by rational measured foliations (being careful
  about Warning~\ref{warn:map-mf-cont}) gives a sequence of simple
  curves (not necessarily connected) so the ratio of extremal lengths
  approaches~$K$.

  The case when there is a conformal embedding, but not a strict
  conformal embedding, can be treated by, for
  instance, adding annuli to the boundary components of~$\Sigma_1$ so
  there is no conformal embedding.
\end{proof}

\begin{corollary}\label{cor:SF-Teich}
  If $\Sigma_1$ and $\Sigma_2$ are closed surfaces and $\phi \co \Sigma_1
  \to \Sigma_2$ is a homeomorphism, then the 
  Teichmüller distance between $\Sigma_1$ and~$\Sigma_2$ is
  $\frac{1}{2} \log \SF[\phi]$, in the sense
  that if $\Sigma_0$ is a fixed base surface and $\psi \co \Sigma_0
  \to \Sigma_1$ is a marking, then the distance between the marked
  surfaces $(\Sigma_1, \psi)$ and $(\Sigma_2, \phi \circ \psi)$ is
  $\frac{1}{2} \log \SF[\phi]$.
\end{corollary}

\begin{proof}
  Immediate from Theorem~\ref{thm:SF-Teich} and the definition
  of Teichmüller distance.
\end{proof}

In this
context, Proposition~\ref{prop:SF-compose} is the triangle
inequality for Teichmüller distance.
  Corollary~\ref{cor:SF-Teich} was
  proved
  by Kerckhoff \cite[Theorem~4]{Kerckhoff80:AsympTeich}. In this case
  connected simple curves suffice.

\begin{remark}
  Despite Warning~\ref{warn:map-mf-cont}, there is a non-continuous
  extension of~$\phi_*$ to measured foliations. Define a pull-back
  function~$\phi^*\co\Curves^-(\Sigma_2) \to
  \powerset\Curves^-(\Sigma_1)$,
  where $\powerset \Curves^-(\Sigma_1)$ is the set of subsets of
  $\Curves^-(\Sigma_1)$,
  by taking all possible intersections
  of~$C$ with the image of~$\Sigma_1$:
  \[
  \phi^*(C) \coloneqq \bigl\{\, [c \cup \phi(\Sigma_1)]
     \bigm| \text{$c$ a representative of $C$}\,\bigr\},
  \]
  where we delete inadmissible components as in the definition
  of~$\phi_*$. Then, for a measured foliation~$F_1 \in \MF^+(\Sigma_1)$,
  we may define $\phi_*(F_2)$ to be the (unique) measured foliation
  $F_2 \in \MF^+(\Sigma_2)$ so that, for all $C_2 \in
  \Curves^-(\Sigma_2)$,
  \[
  i([C_2],F_2) = \inf_{C_1 \in \phi^*(C_2)} i([C_1],F_1).
  \]
  With this definition, the curves in Equation~\eqref{eq:el-surface}
  can be replaced by measured foliations, and the supremum is achieved.
\end{remark}

\begin{definition}
 The embedding $\phi$ is
  \emph{annular} if it extends to an embedding of an annular extension
  $\widehat\Sigma_1$ in~$\Sigma_2$, where $\widehat\Sigma_1$ is
  obtained by attaching an annulus to each boundary component
  of~$\Sigma_1$.
\end{definition}

\begin{theorem}[{Joint with Kahn and Pilgrim \cite[Theorems 1 and 2]{KPT15:EmbeddingEL}}]\label{thm:el-embedding}
  If $(\Sigma_1, \omega_1)$ and $(\Sigma_2, \omega_2)$ are two
  conformal surfaces and $\phi \co \Sigma_1 \hookrightarrow \Sigma_2$
  is a topological embedding, then $\phi$ is homotopic to a conformal
  embedding iff $\SF[\phi] \le 1$.

  Furthermore, the following conditions are equivalent:
  \begin{enumerate}
  \item\label{item:emb-sf-small} $\SF[\phi] < 1$,
  \item\label{item:emb-annular} $\phi$ is homotopic to an annular
    conformal embedding, and
  \item\label{item:emb-strict} $\phi$ is homotopic to a strict
    conformal embedding.
  \end{enumerate}
\end{theorem}

\begin{proof}[Proof sketch]
  The case when $\phi$ is not homotopic to an embedding is implied by
  Theorem~\ref{thm:SF-Teich}. The converse of the first claim is
  essentially Schwarz's Lemma.

  To show that (\ref{item:emb-sf-small}) implies~(\ref{item:emb-annular}),
  pick a quadratic differential $q \in \Quad^+(\Sigma_1)$ that is
  strictly positive on each boundary component.
  Define $\widehat{\Sigma}_1^t$ to be
  $\Sigma_1$ plus an annulus of width~$t$ on each boundary
  component, using the coordinates from~$q$. Elementary estimates
  show that $\SF[\widehat{\Sigma}_1^t \to \Sigma_1]$ (in the natural
  homotopy class) approaches~$1$ as $t$
  approaches~$0$. Then for $t$ small, by Proposition~\ref{prop:SF-compose}
  \[
  \SF[\widehat{\Sigma}_1^t \to \Sigma_2]
    \le \SF[\widehat{\Sigma}_1^t \to \Sigma_1]\cdot \SF[\Sigma_1 \to \Sigma_2]
    < 1,
  \]
  so by Theorem~\ref{thm:SF-Teich} there is a conformal embedding of $\widehat{\Sigma}_1^t$
  in~$\Sigma_2$.

  Clearly (\ref{item:emb-annular}) implies
  (\ref{item:emb-strict}). Finally, if $\phi$ is a strict conformal
  embedding of
  $\Sigma_1$ in~$\Sigma_2$, then by considering test metrics we can
  show that
  \begin{equation}
    \label{eq:area-bound}
    \SF[\phi] \le \sup_{\substack{q \in \Quad^+(\Sigma_2)\\q \ne 0}}
      \frac{\Area_q(\Image(\phi))}{\Area_q(\Sigma_2)}.
  \end{equation}
  Here $\Image(\phi) \subset \Sigma_2$ is the image of $\phi$, which by
  hypothesis misses an open subset of~$\Sigma_2$, and $\Area_q$ is the
  area with respect to the quadratic differential~$q$. Thus for each~$q$,
  $\Area_q(\Image(\phi))/\Area_q(\Sigma_2) < 1$. Since the supremand
  doesn't change as we scale~$q$, we are maximizing over the
  compact set $\PP\Quad^+(\Sigma_2,\omega_2)$ and the supremum is strictly less
  than~$1$.
\end{proof}

\subsection{Behaviour of stretch factor under covers}
\label{sec:sf-cover}

In Section~\ref{sec:asymptotic-SF}, we will also need an
understanding of the behavior of $\SF$ under covers.
Let $\phi \co \Sigma_1 \hookrightarrow \Sigma_2$ be a topological embedding of
  conformal surfaces, let $\wt{\Sigma}_2$ be a finite cover of
  $\Sigma_2$, and let $\wt{\phi} \co \wt{\Sigma}_1 \to
  \wt{\Sigma}_2$ be the pull-back of~$\phi$;
that is, $\wt\Sigma_1$ is defined by the pull-back of $\wt\Sigma_2$
and $\Sigma_1$ via the diagram
\begin{equation}
\begin{tikzpicture}[baseline=($(S1t.base)!0.5!(S1.base)$)]
    \matrix[row sep=0.7cm,column sep=0.8cm] {
      \node (S1t) {$\wt\Sigma_1$}; &
        \node (S2t) {$\wt\Sigma_2$}; \\
      \node (S1) {$\Sigma_1$}; &
        \node (S2) {$\Sigma_2$.}; \\
    };
    \draw[dashed,right hook->] (S1t) to node[above,cdlabel] {\wt\phi} (S2t);
    \draw[dashed,->] (S1t) to (S1);
    \draw[->] (S2t) to (S2);
    \draw[right hook->] (S1) to node[above,cdlabel] {\phi} (S2);
\end{tikzpicture}
\label{eq:pullback}
\end{equation}
When $\wt\Sigma_2$ is connected, $\wt\Sigma_1$ is the minimal cover
of~$\Sigma_1$
with a map to~$\wt\Sigma_2$ making the above diagram commute.

\begin{question}\label{quest:SF-covers}
  How does $\SF[\wt{\phi}]$ compare to $\SF[\phi]$?
\end{question}

It appears that $\SF[\wt{\phi}] \ne \SF[\phi]$ in general,
at least if we allow $\Sigma_2$ to be an
orbifold. However, we can make many partial statements.

\begin{lemma}\label{lem:SF-cover-inc}
  For $\phi$ and $\wt{\phi}$ as above, $\SF[\wt{\phi}] \ge \SF[\phi]$.
\end{lemma}

\begin{proof}
  Follows from the definition of $\SF$ and the good behavior of
  extremal length under covers, Lemma~\ref{lem:el-cover}.
\end{proof}

\begin{lemma}\label{lem:SF-cover-Ioffe}
  For $\phi$ and $\wt{\phi}$ as above, if $\SF[\phi] \ge 1$, then
  $\SF[\wt{\phi}] = \SF[\phi]$.
\end{lemma}

\begin{proof}
  Follows from Theorem~\ref{thm:SF-Teich}.
\end{proof}

\begin{lemma}\label{lem:sf-cover-1}
  For $\phi$ and $\wt{\phi}$ as above, $\SF[\wt\phi] < 1$ iff
  $\SF[\phi] < 1$ and $\SF[\wt\phi] = 1$ iff
  $\SF[\phi] = 1$.
\end{lemma}

\begin{proof}
  If $\SF[\phi] < 1$, by Theorem~\ref{thm:el-embedding} the map~$\phi$
  is homotopic to a
  strict conformal embedding,
  and so $\wt{\phi}$ is too and thus $\SF[\wt\phi] < 1$. The rest follows
  from Lemma~\ref{lem:SF-cover-Ioffe}.
\end{proof}

\begin{definition}
  For $\phi \co \Sigma_1 \to \Sigma_2$ a topological embedding of
  surfaces, define the \emph{lifted
    stretch factor} by
  \[
  \wt\SF[\phi] \coloneqq \lim_{\wt\phi\text{ covers }\phi} \SF[\wt\phi],
  \]
  where the limit runs over increasingly large finite covers of $\phi$
  (determined by a covering of $\Sigma_2$). These covers form a
  directed system, and $\SF$ only increases in a cover
  (Lemma~\ref{lem:SF-cover-inc}) while remaining bounded
  (Lemmas~\ref{lem:SF-cover-Ioffe} and~\ref{lem:sf-cover-1}), so the
  limit exists.
\end{definition}

By definition, if $\wt\phi$ is a covering of $\phi$, then $\wt\SF[\wt\phi]
= \wt\SF[\phi]$.

We will ultimately need to show that certain maps have lifted stretch
factor less than one. We give two arguments, one more elementary
and one proving a stronger result. (See the two proofs of
Proposition~\ref{prop:asf-bound-surf}.)

\begin{definition}
  A homotopy class of topological embeddings $[\phi] \co \Sigma_1 \to
  \Sigma_2$ is \emph{conformally loose} if, for all $y \in \Sigma_2$,
  there is a conformal embedding $\psi \in [\phi]$ so that $y \notin
  \overline{\phi(\Sigma_1)}$.
\end{definition}

  If $\Sigma_2$ is compact
  and $[\phi]$ is conformally loose, then we can find finitely many
  conformal embeddings $\phi_i \in [\phi], i = 1,\dots,n$ so that
  \begin{equation}
  \bigcap_{i=1}^n \overline{\phi(\Sigma_1)} = \emptyset.
  \label{eq:finite-loose}
  \end{equation}

\begin{proposition}\label{prop:sf-loose-bound}
  If $[\phi]$ is conformally loose wih $n$ maps~$\phi_i$ in
  Equation~\eqref{eq:finite-loose}, then $\SF[\phi] \le
  1-1/n$.
\end{proposition}

\begin{proof}[Proof sketch]
  Fix a weighted multi-curve $C$ on $\Sigma_1$, and let $q \in
  \Quad^+(\Sigma_2,\omega_2)$ be the quadratic differential corresponding to
  $\phi(C)$. For at least one $i$, we will have
  \[\frac{\Area_q(\phi_i(\Sigma_1))}{\Area_q(\DD)} \le 1-1/n.\]
  Choosing test metrics as in
  Equation~\eqref{eq:area-bound} shows that $\EL_{\omega_1}[C] \le
  (1-1/n)\EL_{\omega_2}[\phi_*C]$, as desired.
\end{proof}

\begin{corollary}\label{cor:loose-lifted-sf}
  If $[\phi]$ is conformally loose, then $\wt\SF[\phi] < 1$.
\end{corollary}

Corollary~\ref{cor:loose-lifted-sf} is enough to prove
Theorem~\ref{thm:detect-rational}. However, more is true.

\begin{theorem}[{Joint with Kahn and Pilgrim \cite[Theorem 3]{KPT15:EmbeddingEL}}]\label{thm:SF-cover-bound}
  For every strict conformal embedding $\phi \co \Sigma_1
  \hookrightarrow \Sigma_2$,
  \begin{equation*}
    \wt\SF[\phi] < 1.
  \end{equation*}
\end{theorem}

By
Equation~\eqref{eq:area-bound}, Theorem~\ref{thm:SF-cover-bound} is a
consequence of the following proposition.

\begin{proposition}\label{prop:subsurf-area}
  Let $\Sigma_2$ be a Riemann surface, and let $\Sigma_1 \subset
  \Sigma_2$ be a subsurface with compact closure. Then there is a
  constant $K < 1$ so that, for all $0 \ne q \in \Quad(\Sigma_2)$,
  \begin{equation}
    \label{eq:subsurf-area}
    \frac{\Area_q(\Sigma_1)}{\Area_q(\Sigma_2)} < K.
  \end{equation}
  Furthermore, $K$ can be chosen so that
  Equation~\eqref{eq:subsurf-area} holds for all coverings of the pair
  $(\Sigma_1, \Sigma_2)$.
\end{proposition}

Proposition~\ref{prop:subsurf-area}, in turn, depends on the following
local lemma.

\begin{figure}
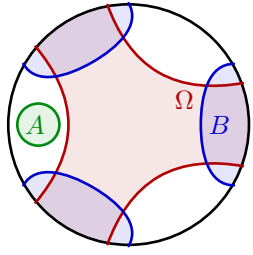

  \[\mfigb{embedding-10}\]
  \caption{The schematic setup of Lemma~\ref{lem:disk-area}.}
  \label{fig:disk-area}
\end{figure}

\begin{lemma}\label{lem:disk-area}
  Let $\Omega \subset \DD$ be an open subset of the disk with an open
  set~$A$ in the complement of $\Omega$, and let $B
  \subset \overline{\DD}$ be a neighborhood of
  $\overline{\Omega} \cap \partial\overline{\DD}$. Then, for every $\epsilon > 0$, there is a
  $\delta > 0$ so that, if $q \in \Quad(\DD)$ is such that
  \begin{equation}
    \label{eq:area-1}
    \frac{\Area_q(\Omega)}{\Area_q(\DD)} > 1-\delta,
  \end{equation}
  then
  \begin{equation}
    \label{eq:area-2}
    \frac{\Area_q(B)}{\Area_q(\DD)} > 1-\epsilon.
  \end{equation}
\end{lemma}

Essentially, Lemma~\ref{lem:disk-area} says that if the area of~$q$ is
concentrating in a subset of the disk, then it is concentrating near
the boundary. See Figure~\ref{fig:disk-area}.

\begin{proof}[Proof sketch for Lemma~\ref{lem:disk-area}]
  If there are no such bounds, there is an $\epsilon$ so that we can find
  a sequence of quadratic differentials $q_n\in\Quad(\DD)$ so that
  \begin{align}
    \Area_{q_n}(\DD) &= 1\\
    \Area_{q_n}(B) &< \epsilon\label{eq:B-small}\\
    \Area_{q_n}(\Omega) &> 1-1/n.\label{eq:omega-large}
  \end{align}
  Consider $\abs{q_n}$ as a measure on $\overline{\DD}$. Since the
  space of measures of unit area on the closed disk is compact in the
  weak topology, after passing to a subsequence we may assume that
  $\abs{q_n}$ converges in the weak topology to some limiting measure
  $\mu$ on~$\overline{\DD}$. Absolute values of holomorphic functions
  are closed in the weak topology, so the restriction of~$\mu$ to the
  open disk can be written $\abs{q_\infty}$ for some holomorphic
  quadratic differential~$q_\infty$. But $\Area_{q_n}(A) < 1/n$, so
  $\Area_{q_\infty}(A) = 0$, so $q_\infty$ is identically~$0$; hence
  $\mu$ is supported on $\bdy\oDD$. Equation~\eqref{eq:omega-large}
  implies that the support of $\mu$ is also contained in
  $\overline{\Omega}$, and hence in
  $\overline{\Omega}\cap\bdy\oDD$. But this contradicts
  Equation~\eqref{eq:B-small}.
\end{proof}

\begin{proof}[Proof sketch for Proposition~\ref{prop:subsurf-area}]
  Divide $\Sigma_2$ by smooth arcs $\alpha_j$ so that the
  complementary regions
  are all disks; let these disks be $U_i$. Consider a quadratic
  differential on $\Sigma_2$ with a very large proportion of its area
  in~$\Sigma_1$. Then most $U_i$ (as weighted by $\Area_q(U_i)$) must have most
  of their $q$-area in $U_i \cap \Sigma_1$. Lemma~\ref{lem:disk-area} then
  says that most of the area of most of the $U_i$ must be in a small
  neighborhood of the seams $\alpha_j$. Arrange the constants so that
  more than half the total area must be in these neighborhoods.

  Now pick an
  alternate set of seams $\beta_k$ with the same
  properties, but so that $\beta_k \cap \alpha_j = \emptyset$ for all
  $j$ and~$k$. By the same argument,
  more than half the total area must
  be concentrated in a neighborhood of the $\beta_k$, a contradiction.

  All of the estimates in this proof only depend on the local geometry
  of the~$U_i$, and thus remain unchanged under taking
  covers.
\end{proof}

\subsection{Stretch factors and embedding energy for graphs}
\label{sec:SF-graphs}

We now turn to the (easier) parallel theory of stretch factors for maps between
graphs. Let $\phi \co \Gamma_1 \to \Gamma_2$
be a continuous map
between two graphs, and suppose that $\Gamma_1$ and~$\Gamma_2$ have
elastic structures $\alpha_1$ and~$\alpha_2$, respectively.

\begin{definition}
  In this setting, the \emph{$\EL$ stretch factor} of~$\phi$ is
  \begin{equation}
    \label{eq:SF-graph}
    \SF_{\EL}[\phi] \coloneqq
      \sup_{\text{$C$ curve on $\Gamma_1$}} \frac{\EL_{\alpha_2}[\phi_* C]}{\EL_{\alpha_1}[C]}.
  \end{equation}
  As for surfaces, this only depends on the homotopy class of~$\phi$.

  Dually, the \emph{Dirichlet stretch factor} of~$\phi$ is
  \begin{equation}
    \label{eq:SF-dir-graph}
    \SF_{\mathrm{Dir}}[\phi] \coloneqq
      \sup_{\substack{\text{$K$ length graph}\\ [f]: \Gamma_2 \to K}} \frac{\Dir[\phi \circ f]}{\Dir[f]},
  \end{equation}
  where the supremum runs over all length graphs $(K,\ell)$ and all
  homotopy classes of maps from $\Gamma_2$ to~$K$.
\end{definition}

\begin{remark}\label{rem:looser}
  As mentioned in the introduction,
  the Dirichlet stretch factor has a natural interpretation in terms
  of rubber-band networks. If $\SF_{\mathrm{Dir}}[\phi] < 1$, then the
  rubber-band network~$\Gamma_1$ is ``looser'' than the rubber-band
  network~$\Gamma_2$: minimal Dirichlet energy in a homotopy class
  with target an arbitrary length graph
  decreases under composition with~$\phi$.
  Theorem~\ref{thm:sf-emb-energy} below implies that this
  is true for arbitrary target geodesic spaces, not just graphs.
\end{remark}

As in the case of surfaces, these stretch factors behave well under
composition (Proposition~\ref{prop:SF-compose}).
Unlike in the case of surfaces, for graphs we have a direct
characterization of the stretch factor.

\begin{definition}\label{def:emb-energy}
  For $\phi \co \Gamma_1 \to \Gamma_2$ a Lipschitz map between elastic
  graphs, the \emph{embedding energy} of~$\phi$ is
  \begin{align}
    \label{eq:emb-energy}
    \Emb(\phi) &\coloneqq \esssup_{y \in \Gamma_2} \Fill_\phi(y)\\
    \Emb[\phi] &\coloneqq \inf_{\psi \in [\phi]} \Emb(\psi),
  \end{align}
  where $\Fill_\phi$ is the filling function defined in
  Equation~\eqref{eq:fill-def}. We take the $L^\infty$ norm of
  $\Fill_\phi$ rather than the $L^1$ norm used for Dirichlet
  energy. This is equivalent to Equation~\eqref{eq:embedding-1}.
\end{definition}

By comparing Equations~\eqref{eq:embedding-1} and~\eqref{eq:lipschitz},
we see that
\begin{equation}
  \label{eq:emb-vs-lipschitz}
  \Emb(\phi) \ge \Lip(\phi).
\end{equation}

\begin{theorem}\label{thm:sf-emb-energy}
  For $\phi \co \Gamma_1 \to \Gamma_2$ a continuous map between
  elastic graphs,
  \[
  \SF_{\EL}[\phi] = \SF_{\mathrm{Dir}}[\phi] = \Emb[\phi].
  \]
\end{theorem}

Theorem~\ref{thm:sf-emb-energy} should be thought of as analogous to
Theorem~\ref{thm:SF-Teich}, although it applies in all cases, not just
when there fails to be a conformal embedding.
It is also analogous to the relation between Lipschitz energy as
maximum derivative (Equation~\eqref{eq:lipschitz}) and as ratio of
curve lengths (Equation~\eqref{eq:lip-ratio}).

\begin{proof}[Proof sketch]
  For $C$ a curve on $\Gamma_1$ and $\phi\co \Gamma_1 \to \Gamma_2$ a
  Lipschitz map, by pushing forward the scaling function $\rho$ in
  Equation~\eqref{eq:el-graph} we can see that
  \begin{equation}
    \label{eq:emb-el}
    \EL[\phi_* C] \le \EL[C] \cdot \Emb(\phi).
  \end{equation}
  This immediately implies that $\SF_{\EL}[\phi] \le \Emb[\phi]$.

  Similarly, the definition of the Dirichlet energy as the $L^1$ norm
  of the filling function and the embedding energy as the $L^\infty$
  norm of the filling function makes it clear that, for $f \co
  \Gamma_2 \to K$ any map to a length graph,
  \begin{equation}
    \label{eq:emb-dir}
    \Dir(f \circ \phi) \le \Dir(f) \Emb(\phi),
  \end{equation}
  which implies that $\SF_{\mathrm{Dir}}[\phi] \le \Emb[\phi]$.

  To prove the opposite inequalities, we find and analyze a
  representative for the homotopy class $[\phi]$ for which Equations
  \eqref{eq:emb-el} and~\eqref{eq:emb-dir} are tight. We defer the
  description of these nice representatives (analogous to harmonic
  maps, Definition~\ref{def:harmonic-graphs}) until
  Section~\ref{sec:lambda-filling}.
\end{proof}

Motivated by Theorem~\ref{thm:sf-emb-energy}, we make the
following definition.

\begin{definition}
  For $(\Gamma_1,\alpha_1)$ and $(\Gamma_2,\alpha_2)$ two elastic
  strip graphs, a map $\phi \co \Gamma_1 \to \Gamma_2$ is
  \emph{loosening} if $\Emb(\phi) \le 1$. The map is \emph{strictly
    loosening} if $\Emb(\phi) < 1$. We likewise say that a homotopy class
  $[\phi]$ is {(strictly) loosening} if there is a (strictly)
  loosening map
  in~$[\phi]$.
\end{definition}

In contrast with the case for surfaces (Section~\ref{sec:sf-cover}),
it is easy to see how the stretch factor for graphs behaves under
covers.

\begin{proposition}\label{prop:SF-cover-graph}
  Let $\phi \co \Gamma_1 \to \Gamma_2$ be a map of elastic graphs,
  $\wt{\Gamma}_2$ be a cover of~$\Gamma_2$, and
  $\wt{\phi} \co \wt{\Gamma}_1 \to \wt{\Gamma}_2$
  be the pull-back map, defined as in Equation~\eqref{eq:pullback}. Then
  \[
  \SF[\wt{\phi}] = \SF[\phi].
  \]
\end{proposition}

\begin{proof}
  The definition of $\SF[\phi]$ by the maximal expansion of extremal length
  or Dirichlet energy shows that $\SF[\wt \phi] \ge \SF[\phi]$ (as in
  Lemma~\ref{lem:SF-cover-inc}). The definition of $\Emb[\phi]$ by the
  infimum over all representatives shows that $\Emb[\wt \phi] \le
  \Emb[\phi]$. Theorem~\ref{thm:sf-emb-energy} shows that we have equalities.
\end{proof}

\subsection{Relating graphs and surfaces}
\label{sec:sf-graph-surface}

\begin{figure}
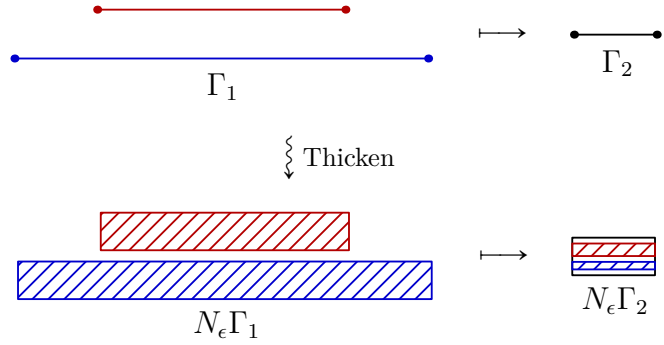

  \begin{gather*}
  \underset{\textstyle\strut\Gamma_1}{\mfigb{surface-21}} \quad\longmapsto\quad
  \underset{\textstyle\strut\Gamma_2}{\mfigb{surface-20}}\\[5pt]
  \mathcenter{\rotatebox{-90}{$\rightsquigarrow$}}
    \footnotesize{\,\,\text{Thicken}}\\[5pt]
  \underset{\textstyle\strut N_\epsilon\Gamma_1}{\mfigb{surface-23}}
    \quad\longmapsto\quad
  \underset{\textstyle\strut N_\epsilon\Gamma_2}{\mfigb{surface-22}}
  \end{gather*}
  \caption{Motivation for the definition of embedding energy. Two
    edges of $\Gamma_1$ mapping to one edge of $\Gamma_2$ get
    thickened up to two rectangles
    mapping to one rectangle. For the result to be a conformal
    embedding, we need the total height of the image to be less than
    the total height available in the range.}
  \label{fig:fit-rectangles}
\end{figure}

As motivation for the somewhat strange definition of embedding energy,
consider two elastic graphs
$(\Gamma_1,\alpha_1)$ and~$(\Gamma_2,\alpha_2)$ and a conformal
embedding of the thickenings $N_t \Gamma_1 \hookrightarrow
N_t\Gamma_2$. Suppose that this conformal embedding is ``close'' in
some sense to a graph map $\phi\co\Gamma_1 \to \Gamma_2$.
If we look away from the thickenings of the vertices (of
both graphs), we see locally a map from some edges of $\Gamma_1$ into
a single edge of~$\Gamma_2$. The total height of images of the
thickened edges of~$\Gamma_1$ must be less than or equal to the total
available height in the thickened edge of~$\Gamma_2$. Near a point $x
\in \Gamma_1$, the image is stretched horizontally by a factor of
$\abs{\phi'(x)}$; since we are considering a conformal embedding, the
image must also be stretched vertically by the same factor. Thus the
image of this portion of the edge takes up a height of
$t\abs{\phi'(x)}$. This argument suggests that if there is a conformal
embedding close to $\phi$, we must have, for each $y \in \Gamma_2$ not
near a vertex or the image of a vertex,
\begin{equation*}
  \sum_{\phi(x) = y} t\abs{\phi'(x)} \le t\qquad\text{or, equivalently,}\qquad
  \Fill_\phi(y) \le 1.
\end{equation*}
See Figure~\ref{fig:fit-rectangles}.

We will not attempt to make the above heuristic argument precise.
Instead, we get a precise statement another way.

\begin{proposition}\label{prop:sf-graph-surface}
  Let $\Gamma_1$ and $\Gamma_2$ be trivalent
  elastic ribbon graphs, and let $m = \min_e \alpha(e)$ be the
  smallest weight of any edge in either $\Gamma_1$ or~$\Gamma_2$. Let
  $\phi \co \Gamma_1 \to \Gamma_2$ be a map that extends to a topological
  embedding $N_t\phi \co N_t\Gamma_1 \to N_t\Gamma_2$.
  Then, for $t < m/2$, we have
  \[
  \SF[\phi]/(1+8t/m) \le
    \SF[N_t\phi] \le
    \SF[\phi]\cdot (1 + 8t/m).
  \]
\end{proposition}

\begin{proof}
  Immediate from Proposition~\ref{prop:el-graphs-surfaces}.
\end{proof}

\begin{theorem}\label{thm:embedding-graph-surface}
  Let $\Gamma_1$ and $\Gamma_2$ be elastic ribbon
  graphs, and let $\phi\co
  \Gamma_1 \to \Gamma_2$ be a map of graphs
  that extends to a topological embedding $N\phi\co N\Gamma_1 \to
  N\Gamma_2$. If $\Emb[\phi] < 1$, then for all $t$
  sufficiently small, $N_t\Gamma_1$ conformally embeds in $N_t\Gamma_2$
  in $[N\phi]$. Conversely, if $N_t\Gamma_1$ conformally embeds
  in $N_t\Gamma_2$ in $[N\phi]$ for all sufficiently small~$t$, then
  $\Emb[\phi] \le 1$.
\end{theorem}

\begin{proof}
  Follows from Proposition~\ref{prop:sf-graph-surface} and
  Theorem~\ref{thm:el-embedding}. If the graphs $\Gamma_1$ and
  $\Gamma_2$ are not trivalent, we
  may approximate them by trivalent graphs at the cost of a small
  change in the stretch factor.
\end{proof}

\subsection{Minimizing embedding energy: $\lambda$-filling maps}
\label{sec:lambda-filling}

We say a little more about the proof of
Theorem~\ref{thm:sf-emb-energy}, and in particular characterize the optimal
maps. Recall from Section~\ref{sec:graphs} that a strip graph has both
lengths and widths.

\begin{definition}
  For $\lambda > 0$, a map $\phi \co S_1 \to S_2$ between
  strip graphs is
  \emph{$\lambda$-filling} if
  \begin{enumerate}
  \item $\phi$ is length-preserving: it does not
    backtrack, and $\abs{\phi'(x)} = 1$ for almost all~$x$, where we take
    derivative with respect to the length metrics, and
  \item $\phi$ scales widths by~$\lambda$: for almost all $y \in S_2$,
    \[
    \sum_{\phi(x)=y} w(x) = \lambda \cdot w(y).
    \]
  \end{enumerate}
\end{definition}

Note that for maps between strip graphs, we have to be careful
whether we differentiate with respect to the length metric or with
respect to the coordinates from the elastic weights, as in
Warning~\ref{warn:elastic-vs-length}.

\begin{lemma}
  If $\phi \co S_1 \to S_2$ is a $\lambda$-filling map
  between strip graphs, then the underlying map
  $\Gamma(\phi) \co \Gamma(S_1) \to \Gamma(S_2)$
  between elastic graphs has embedding energy $\lambda$. This is
  minimal in the homotopy class.
\end{lemma}

\begin{proof}[Proof sketch]
  The conditions on~$\phi$ imply that $\Fill_{\Gamma(\phi)}=\lambda$
  everywhere on $\Gamma(S_2)$. (Recall that ``length'' on
  the elastic graph is $\text{length}/\text{width}$ in terms of the strip
  graph.)

  Now consider the length graph~$K(S_2)$. The map $\Gamma(S_2)
  \to K(S_2)$ has Dirichlet energy $\Area(S_2)$, almost by definition. The
  composite map $\Gamma(S_1) \to \Gamma(S_2) \to K(S_2)$ has Dirichlet
  energy $\lambda\cdot \Area(S_2)$. It follows that
  $\SF_{\Dir}[\Gamma(\phi)] \ge \lambda$. The easy direction of
  Theorem~\ref{thm:sf-emb-energy} then implies that
  $\Emb[\Gamma(\phi)] = \lambda$.
\end{proof}

We
say that a map~$\phi$ between elastic graphs is \emph{$\lambda$-filling} if
there are compatible strip structures that make~$\phi$
$\lambda$-filling.
It is not true that every homotopy class of maps between elastic
graphs has a $\lambda$-filling representative. However, we can make it
$\lambda$-filling on a subgraph.

\begin{definition}
  A map $\phi \co S_1 \to S_2$ between strip graphs is
  \emph{partially $\lambda$-filling} if there are non-empty subgraphs $T_1$ of
  $S_1$ and $T_2$ of $S_2$ so that
  \begin{enumerate}
  \item $\phi(T_1) = T_2$ and $\phi^{-1}(T_2) = T_1$;
  \item the restriction of $\phi$ to a map $T_1 \to T_2$ is
    $\lambda$-filling;
  \item $\phi$ is everywhere length-preserving; and
  \item outside of $T_1$ and $T_2$, the map $\phi$ scales widths by less
    than~$\lambda$.
  \end{enumerate}
\end{definition}

\begin{proposition}\label{prop:partially-lambda}
  Every homotopy class of maps between elastic graphs has a partially
  $\lambda$-filling representative.
\end{proposition}

Theorem~\ref{thm:sf-emb-energy} follows quickly from
Proposition~\ref{prop:partially-lambda}.

\subsection{Computing embedding energy}
\label{sec:compute-emb}
Given two
elastic graphs $\Gamma_1$ and~$\Gamma_2$ and a homotopy class $[\phi]
\co \Gamma_1 \to \Gamma_2$ of maps between them, how can we concretely
find the partially $\lambda$-filling
representative $\psi \in [\phi]$ guaranteed by
Proposition~\ref{prop:partially-lambda} and thus compute $\Emb[\phi]$?
The following iteration
appears to converge.
\begin{enumerate}
\setcounter{enumi}{-1}
\item Pick a set of widths $v_0 \in \Wid(\Gamma_2)$, i.e., a width for
  each edge of $\Gamma_2$ satisfying the triangle inequality.
\item\label{item:step-K} Take the metric graph $K_0$ to be $\Gamma_2$ with edge $e$
  assigned length $\alpha(e) v_0(e)$. Note that the evident map $f_0 \co
  \Gamma_2\to K_0$ is harmonic.
\item Find a harmonic representative~$g_0$ of the composite map $\Gamma_1
  \overset{\phi}{\to} \Gamma_2 \overset{f}{\to} K_0$. Our first
  approximation to~$\psi$ is $\psi_0 = (f_0)^{-1} \circ g_0$.
\item Compute half the tension of~$g_0$ (i.e., $\abs{g_0'(x)}$) in each
  edge of~$\Gamma_1$, getting $w_0 \in \Wid(\Gamma_1)$.
\item Push forward $w_0$ to a function~$v_1 \in \Wid(\Gamma_2)$ by
  setting, for $e \in \Edges(\Gamma_2)$ and $y \in e$,
  \[
  v_1(y) = \sum_{x \in \psi_0^{-1}(y)} w_0(x).
  \]
  This is independent of the choice of~$y$ since $g$ is harmonic. Then
  return to Step~\eqref{item:step-K}, using $v_1$ instead of~$v_0$.
\end{enumerate}
Schematically, we are iterating around the following cycle:
\[
\begin{tikzpicture}[baseline]
  \matrix[row sep=20pt, column sep=2em]{
    \node (L2) {$\Len(\Gamma_2)$};
      & \node (W2) {$\Wid(\Gamma_2)$};\\
    \node (L1) {$\Len(\Gamma_1)$};
      & \node (W1) {$\Wid(\Gamma_1).$};\\
  };
  \draw [->] (L2) -- node[left]{\small (2)} (L1);
  \draw [->] (L1) -- node[above]{\small (3)} (W1);
  \draw [->] (W1) -- node[right]{\small (4)} (W2);
  \draw [->] (W2) -- node[above]{\small (1)} (L2);
\end{tikzpicture}
\]
After each iteration, we can compute the embedding energy:
\[
\Emb(\psi_i) = \max_{e \in \Edges(\Gamma_2)} \frac{v_{i+1}(e)}{v_i(e)}.
\]

\begin{conjecture}
  The algorithm above converges to a map with lowest embedding energy.
\end{conjecture}
Once the combinatorics of the graph map have settled down, this
maps in this iteration become linear, and the algorithm reduces to
finding the largest eigenvector of a linear system by iteration.
In practice, the algorithm appears to converge rapidly.

Note the relation to tightness in
Equation~\eqref{eq:norm-submultiplicative}: we simultaneously find a
representative $\psi \in [\phi]$, a map $f$ from $\Gamma_2$ to a
metric graph~$K$, and a set of widths~$w$ on~$\Gamma_1$, all
multiplicative on the nose:
\[
\Emb(\psi) = \frac{\Dir[f \circ \psi]}{\Dir[f]} = \frac{\EL[\psi_* w]}{\EL[w]},
\]
where we use a natural extension of $\EL$ to a function on $\Wid(\Gamma)$.

\subsection{Extensions and questions}
\label{sec:sf-extensions}

\begin{conjecture}\label{conj:cover-loose}
  For any strict conformal embedding $\phi \co \Sigma_1 \to \Sigma_2$,
  there is a cover $\wt \phi \co \wt\Sigma_1 \to \wt\Sigma_2$ so that
  $[\wt\phi]$ is conformally loose.
\end{conjecture}

\begin{question}
  What happens if we vary the definition of $\SF$ for surfaces? For
  instance, we restricted to simple curves in
  Definition~\ref{def:sf-surface}. What happens if we drop that
  restriction, and look at general curves? What if we look at the
  expansion factor of Dirichlet energy for maps to graphs instead?
  What if we look at
  width graphs, as in Problem~\ref{prob:el-emb-graphs}?
\end{question}

\begin{problem}\label{prob:sf-surface}
  Give a direct expression for (some version of) the stretch factor
  $\SF[\phi]$ or lifted stretch factor $\wt\SF[\phi]$ of
  a conformal surface embedding, analogous to Theorem~\ref{thm:SF-Teich} for
  the case when $\phi$ is not homotopic to a conformal embedding or
  to Theorem~\ref{thm:sf-emb-energy} for the case of graphs.
\end{problem}

\begin{problem}\label{prob:sf-fill-graph}
  In Definition~\ref{def:el-graph-surf}, we defined a notion of
  extremal length starting from an elastic graph $\Gamma$ embedded as
  a filling subset of a surface. Give a direct expression for the
  stretch factor (maximal ratio of extremal lengths) between two
  graphs $\Gamma_1$ and~$\Gamma_2$ with filling embeddings in the same
  surface~$\Sigma$. (Theorem~\ref{thm:sf-emb-energy} handles the case
  when $\Gamma_2$ is a spine.)
\end{problem}

Since filling graphs can be used to approximate surfaces, a solution
to Problem~\ref{prob:sf-fill-graph} would presumably be helpful in
answering Problem~\ref{prob:sf-surface}.

\begin{problem}
  Extend Problem~\ref{prob:sf-fill-graph} to the more general setting
  of a group~$G$ and surjective maps $\pi_1(\Gamma_1) \to G$ and
  $\pi_1(\Gamma_2) \to G$, as in Definition~\ref{def:el-group}.
\end{problem}


\section{Dynamics}
\label{sec:dynamics}

\subsection{Iterating covers}
\label{sec:iterating-covers}
We finally turn to the \emph{dynamical} picture: what happens when we
iterate a map from a conformal surface or elastic graph to itself? Our
setting is not quite the usual dynamical picture: we are
``iterating'' virtual endomorphisms, which are not maps
from a space $X$ to itself, but maps from a cover $\wt X$ to $X$.

\begin{definition}\label{def:iterated-covers}
  Let $X_0$ and $X_1$ be topological spaces, $\pi \co X_1 \to X_0$ be a covering
  map of degree~$d$,
  and $\phi \co X_1 \to X_0$ be a continuous map. We call this data a
  \emph{virtual endomorphism} of~$X_0$, also called a \emph{topological
  automaton} by Nekrashevych
\cite{Nekrashevych05:SelfSimilar,Nekrashevych13:CombModel}.
  Define $X_k$ to be the $k$-fold
  product of $X_1$ with itself over $X_0$ using the two maps $\pi$
  and~$\phi$; e.g., $X_3$ is the pullback of the diagram below.
  \[
  \begin{tikzpicture}[node distance=1.5cm]
    \node (X00) at (0,0) {$X_0$};
    \node (X20) at (2,0) {$X_0$};
    \node (X40) at (4,0) {$X_0$};
    \node (X60) at (6,0) {$X_0$};

    \node (X11) at (1,1) {$X_1$};
    \node (X31) at (3,1) {$X_1$};
    \node (X51) at (5,1) {$X_1$};


    \node (X33) at (3,3) {$X_3$};

    \draw[->] (X11) to node[above left=-1mm,cdlabel]{\pi} (X00);
    \draw[->] (X31) to node[above left=-1mm,cdlabel]{\pi} (X20);
    \draw[->] (X51) to node[above left=-1mm,cdlabel]{\pi} (X40);

    \draw[->] (X11) to node[above right=-1mm,cdlabel]{\phi} (X20);
    \draw[->] (X31) to node[above right=-1mm,cdlabel]{\phi} (X40);
    \draw[->] (X51) to node[above right=-1mm,cdlabel]{\phi} (X60);


    \draw[->,dashed] (X33) to (X11);
    \draw[->,dashed] (X33) to (X31);
    \draw[->,dashed] (X33) to (X51);

    \draw[->,dashed,bend right] (X33) to node[above left=-1mm,cdlabel]{\pi_3} (X00);
    \draw[->,dashed,bend left] (X33) to node[above right=-1mm,cdlabel]{\phi_3} (X60);
  \end{tikzpicture}
  \]
  Concretely, define
  \begin{equation*}
    X_k \coloneqq \bigl\{\,(x_1,\dots,x_k) \in (X_1)^k \bigm| \phi(x_i) =
       \pi(x_{i+1}) \,\bigr\}
  \end{equation*}
  $X_k$ comes with two natural maps to~$X_0$:
  \begin{itemize}
  \item The map $\pi_k$ is the map to the leftmost copy
    of~$X_0$:
    \[
    \pi_k(x_1,\dots,x_k) \coloneqq \pi(x_1).
    \]
    It is a covering map of degree~$d^k$.
  \item The map $\phi_k$ is the map to the rightmost factor
    of~$X_0$:
    \[
    \phi_k(x_1,\dots,x_k) \coloneqq \phi(x_k).
    \]
    It is a composition of lifts of
    $\phi$ to various covers of~$X_0$.
  \end{itemize}
  We will also use the map $\phi$ to represent the entire virtual
  endomorphism.
\end{definition}

This construction makes sense even if $\pi$ is not a
covering map; in this generality, we are composing topological
correspondences.

\subsection{Asymptotic stretch factors}
\label{sec:asymptotic-SF}

Now consider the case that $X$ is either a conformal surface
$(\Sigma,\omega)$ or an elastic graph $(\Gamma, \alpha)$, with a
virtual endomorphism $\pi,\phi \co
X_1 \rightrightarrows X$ as above. If $X$ is a
conformal surface, suppose also that $\phi$ is a topological
embedding. Note that $X_k$ inherits the structure of a conformal
surface or elastic graph from~$X$ via the covering map $\pi_k$.

\begin{definition}
  The \emph{asymptotic stretch factor} of the virtual endomorphism is
  \begin{equation}
    \label{eq:asf}
    \ASF[\phi] \coloneq \lim_{n \to \infty} \SF[\phi_n]^{1/n}.
  \end{equation}
\end{definition}

\begin{lemma}\label{lem:asf-exists}
  The limit in Equation~\eqref{eq:asf} exists.
\end{lemma}

\begin{proof}[Proof sketch]
  For graphs, sub-multiplicativity of $\SF$
  (Proposition~\ref{prop:SF-compose}) and good behavior under covers
  (Proposition~\ref{prop:SF-cover-graph}) show that
  \[
  \SF[\phi_{k + l}] \le \SF[\phi_{k}] \cdot \SF[\phi_{l}],
  \]
  which implies that the sequence $\SF[\phi_{k}]^{1/k}$
  converges. (It is close to a decreasing sequence.)
  For surfaces,
  $\SF$ does not behave as well under covers, but we still have
  that if $\wt \phi$ is a cover of $\phi$, then $\SF[\phi] \le \SF[\wt
  \phi] \le \max(1,\SF[\phi])$ by Lemmas~\ref{lem:SF-cover-inc},
  \ref{lem:SF-cover-Ioffe}, and~\ref{lem:sf-cover-1}; this is enough
  to show convergence.
\end{proof}

\begin{lemma}\label{lem:asf-indep}
  The asymptotic stretch factor $\ASF$ is independent of the conformal
  structure~$\omega$ or elastic structure~$\alpha$ used to define
  it. More generally, if $f \co \Gamma \to \Gamma'$ is a homotopy
  equivalence of elastic graphs with homotopy inverse $g \co \Gamma'
  \to \Gamma$, then
  \begin{equation}\label{eq:asf-conj}
  \ASF\bigl[f \circ \phi \circ g_1\bigr] = \ASF[\phi],
  \end{equation}
  where $g_1\co \Gamma'_1 \to
  \Gamma_1$ is the lift of~$g$ to the cover.
\end{lemma}

\begin{proof}[Proof sketch]
  Consider Equation~\eqref{eq:asf-conj} for graphs and let $\phi' = f
  \circ \phi \circ g_1$. We have $[\phi'_{n}] = [f \circ \phi_{n} \circ
  g_n]$. Therefore
  \[
  \SF[\phi'_{n}] \le \bigl(\SF[f]\SF[g_n]\bigr) \SF[\phi_{n}],
  \]
  with a similar inequality the other way. Since $\SF[g_n] = \SF[g]$
  (Proposition~\ref{prop:SF-cover-graph}),
  as $n \to \infty$ the contribution of the factor
  $\SF[f]\SF[g_n]$ to the limit goes to~$1$, so
  the asymptotic stretch factors are equal. The case of surfaces is similar.
\end{proof}

Thus we may speak about the asymptotic stretch factor of a virtual
endomorphism $\pi,\phi\co X_1 \rightrightarrows X$ where $X$ is a
topological graph or a surface,
without reference to the conformal or elastic structure. When the
covering map is trivial (i.e., the virtual endomorphism is an ordinary
endomorphism), this recovers W.~Thurston's theory of
pseudo-Anosov maps.

\begin{proposition}\label{prop:asf-teich}
  If $\phi \co \Sigma\righttoleftarrow$ is a pseudo\hyp Anosov
  self\hyp homeomorphism of a
  surface~$\Sigma$ (possibly with boundary), then $\ASF[\phi]$ is the
  pseudo\hyp Anosov
  constant of~$\phi$, i.e., the exponential of the translation
  distance of the induced map on Teichmüller space.
\end{proposition}

\begin{proof}[Proof sketch]
  Follows from Theorem~\ref{thm:SF-Teich}.
\end{proof}

\begin{proposition}\label{prop:asf-graph-surface}
  Let $\Gamma$ be a ribbon graph, and let $\pi,\phi \co \Gamma_1
  \rightrightarrows \Gamma$ be a
  virtual endomorphism of~$\Gamma$ so that $\phi$ extends to a topological
  embedding
  $N\phi \co N\Gamma_1 \hookrightarrow N\Gamma$. Then
  \[
  \ASF[\phi] = \ASF[N\phi].
  \]
\end{proposition}

\begin{proof}[Proof sketch]
  Proposition~\ref{prop:sf-graph-surface} relates stretch factors on
  graphs and on surfaces. The errors in the estimates disappear in the
  limit defining $\ASF$, as in Lemma~\ref{lem:asf-indep}.
\end{proof}

\begin{proposition}\label{prop:asf-bound-graph}
  Let $\Gamma$ be a graph and $\pi,\phi\co\Gamma_1 \rightrightarrows \Gamma$ be a
  virtual endomorphism of~$\Gamma$. Then $\ASF[\phi] \le
  \SF[\phi]$.
\end{proposition}

\begin{proof}
  Immediate from Proposition~\ref{prop:SF-compose},
  Proposition~\ref{prop:SF-cover-graph}, and the definition of
  $\ASF$.
\end{proof}

\begin{proposition}\label{prop:asf-bound-surf}
  Let $\Sigma$ be a conformal surface and $\pi,\phi \co \Sigma_1
  \rightrightarrows \Sigma$ be a
  virtual endomorphism of~$\Sigma$ with $\phi$ a conformal annular
  embedding. Then
  $\ASF[\phi] < 1$.
\end{proposition}

\begin{proof}[Proof, version 1]
  By Theorem~\ref{thm:SF-cover-bound}, $\wt\SF[\phi] < 1$. Then we have
  \[
  \SF[\phi_n] = \SF[\phi_{n-1} \circ \wt\phi] \le \SF[\phi_{n-1}]\wt\SF[\phi]
  \]
  (where $\wt\phi \co X_n \to X_{n-1}$ is a cover of~$\phi$), so
  \begin{align*}
    \SF[\phi_n] &\le \bigl(\wt\SF[\phi]\bigr)^n\\
    \ASF[\phi] &\le \wt\SF[\phi]. \qedhere
  \end{align*}
\end{proof}

\begin{proof}[Proof, version 2]
  Here we avoid Theorem~\ref{thm:SF-cover-bound}.

  Let $J(\phi)$ be the
  \emph{Julia set} of $\phi$: the intersection of the images of
  $\phi_n(\Sigma_n)$. This set has measure~$0$. Suppose for simplicity
  that $\Sigma$ is connected; an Euler characteristic argument shows that
  $\Sigma$ is planar. Pick a conformal embedding of~$\Sigma$ in
  $\CC\PP^1$. For each $x \in
  J(\phi)$ and $n$ sufficiently large, $\phi_n(\Sigma_n)$ can be
  translated in $\Sigma\subset\CC\PP^1$ to a map that misses~$x$. By
  compactness of $J(\phi)$, for sufficiently
  large~$n$ the homotopy class $[\phi_n]$ is conformally loose, so
  $\wt\SF[\phi_n] < 1$ by Proposition~\ref{prop:sf-loose-bound}. This
  then implies $\ASF[\phi] < 1$ as above.
\end{proof}

\begin{proposition}\label{prop:asf-to-sf}
  For $\phi$ a virtual endomorphism of either an elastic graph or a
  conformal surface, if
  $\ASF[\phi] < 1$, then for all sufficiently large~$n$,
  $\SF[\phi_n] < 1$.
\end{proposition}

Loosely speaking, Proposition~\ref{prop:asf-bound-surf} says that
$\ASF[\phi]$ detects
conformal embeddings in covers. Since $\ASF$ is also independent of the
elastic or conformal structure used to define the stretch factor, this
becomes a
useful tool for studying rational maps.

Example~\ref{examp:asf-bad} gives a case where the stretch factor
necessarily drops under iteration, so $\ASF[\phi] < \SF[\phi]$.
  
\subsection{Comparison to Lipschitz expansion}
\label{sec:asf-lipschitz}

There is another natural dynamical invariant of a virtual endomorphism
of a graph: its (asymptotic) Lipschitz expansion, which we now compare
with the asymptotic stretch
factor.

\begin{definition}
  For $\Gamma$ a metric graph and $\pi,\phi \co \Gamma_1 \rightrightarrows \Gamma$ a
  virtual endomorphism of~$\Gamma$, the Lipschitz energy of~$\phi$ was
  defined in
  Equation~\eqref{eq:lipschitz}. The \emph{asymptotic expansion}
  of~$\phi$ is
  \[
  \ALip[\phi] \coloneqq \lim_{n \to \infty} \Lip[\phi_{n}]^{1/n}.
  \]
  This is independent of the metric
  graph used to define it (up to homotopy equivalence), by the
  argument in Lemma~\ref{lem:asf-indep}, .
\end{definition}

A \emph{metric train track representative} for a free group automorphism
$\psi \co F_n \to F_n$ is an endomorphism $\phi$ of a graph $\Gamma$
inducing $\psi$ on~$\pi_1$ so that $\Lip[\phi^n] = \Lip(\phi)^n$,
i.e., so that $\ALip[\phi] = \Lip[\phi] = \Lip(\phi)$.
 Bestvina and
  Handel's theory of
  train tracks for free group automorphisms
  \cite{BH92:TrainTracks,Bestvina11:Bers}
  shows that most free group automorphisms have metric train track
  representatives. (More precisely, any irreducible automorphism has a
  metric train track representative.)

Given the similarity of the definitions, one might suspect these two
quantities are related. Indeed, they are essentially the same in the case of
endomorphisms that are not virtual.

\begin{proposition}\label{prop:lip-emb-1}
  For $\phi \co \Gamma \to \Gamma$ an endomorphism of a graph that
  induces an automorphism on $\pi_1$, we have
  \[
  \ASF[\phi] = \ALip[\phi]^2.
  \]
\end{proposition}

\begin{proof}[Proof sketch]
  Suppose for simplicity that $\Gamma$ is trivalent and
  has $k$ edges each of length~$1$. Then, comparing
  Equations~\eqref{eq:len-graph} and~\eqref{eq:el-graph-2},
  we see that, for any curve~$C$,
  \[
  \frac{1}{k} \ell[C]^2 \le \EL[C] \le \ell[C]^2.
  \]
  We deduce, from Equation~\eqref{eq:lip-ratio} and the definition of
  stretch factor, that
  \begin{equation}\label{eq:lip-sf}
  \frac{1}{k} \Lip[\phi]^2 \le \SF[\phi] \le k\Lip[\phi]^2.
  \end{equation}
  The error factors in these inequalities go away in the limits
  defining $\ASF$ and $\ALip$.
\end{proof}

\begin{remark}
  Proposition~\ref{prop:lip-emb-1} may be compared to the fact that
  the absolute value of the largest eigenvalue of a linear operator on
  a finite-dimensional space equals the asymptotic growth rate of the
  norm of any vector, independent of the choice of norm. In this
  analogy, $\EL$ is like the (square of the) $L^2$ norm while
  $\ell$ is like the $L^1$ norm.
\end{remark}

For a \emph{virtual} endomorphism $\pi,\phi\co \Gamma_1
\rightrightarrows \Gamma$, the
quantities $\ASF$ and $\ALip$ are in general different. By
Equation~\eqref{eq:emb-vs-lipschitz}, we have the rather weak
inequality
\begin{align}
  \ALip[\phi] &\le \ASF[\phi].
  \label{eq:ASF-ALip}
\end{align}

See Section~\ref{sec:branched-lipschitz} for a comparison of how these
two quantities
relate to dynamics on~$S^2$.

\begin{remark}
  The reason it is impossible to find train-track representatives in
  cases like
  Example~\ref{examp:asf-bad} is that the same vertex in~$\Gamma$
  appears multiple times in $\Gamma_1$, and it is impossible to
  arrange for all of these vertices to simultaneously obey the
  train-track requirements.
\end{remark}

\subsection{Extensions and questions}
\label{sec:asf-extensions-questions}

There are many questions raised by this theory. First of all, in order
to cover the general case of rational maps we need to extend the
definitions to orbifold fundamental groups.

\begin{problem}
  Extend the definition of $\ASF$ to automorphisms of pointed
  graphs, or to virtual endomorphisms of orbifolds.
\end{problem}

\begin{problem}
  More generally, define and study $\ASF$ for a virtual
  endomorphism of an arbitrary group.
\end{problem}

We can also ask about properties of the result.

\begin{question}\label{quest:asf-algebraic}
  Is $\ASF[\phi]$ algebraic?
\end{question}

\begin{problem}\label{prob:find-ASF}
  Give an algorithm for computing $\ASF[\phi]$, either approximately
  or exactly.
\end{problem}


\section{Rational maps}
\label{sec:rational-maps}

Finally, we come to the original goal of this work, the study of
branched self-covers of the sphere and when they are equivalent to
rational maps.

\begin{definition}\label{def:branched-self-cover}
  Fix a finite set~$P$ of points in a sphere~$S^2$. A \emph{branched
    self-cover} $f \co (S^2,P)\righttoleftarrow$ is a map so that
  $f(P)\subset P$ and so that $f$ is a covering map when restricted to
  $S^2\setminus f^{-1}(P) \to S^2\setminus P$. Two branched self-covers
  are \emph{equivalent} if they are related by homotopy relative
  to~$P$ and by conjugacy. (Without conjugacy, we could not change the
  set~$P$.) A branched self-cover is \emph{rational} if it is
  equivalent to a rational map on $\CC\PP^1$, which is necessarily
  post-critically finite. ($P$ contains the post-critical set, but may
  be larger.)
\end{definition}

We will assume there is a branch point in each cycle of~$P$. For
rational maps, this implies $f$ is hyperbolic.

\subsection{Characterizing rational maps conformally}
\label{sec:rational-conformal}

We first give a characterization of rational maps in terms of
surfaces.

\begin{definition}

  A \emph{surface spine} for $S^2\setminus P$ is a surface~$\Sigma$
  (necessarily planar) together with a topological embedding $i \co
  \Sigma \hookrightarrow S^2\setminus P$ so that complement of the
  image consists of one punctured disk for each point in~$P$, i.e., so
  that $S^2\setminus P$ deformation retracts onto~$i(\Sigma)$.
\end{definition}

\begin{definition}
  Let $f\co (S^2,P)\righttoleftarrow$ be a branched self-cover and
  $\Sigma$ be a surface equipped with an embedding $i \co \Sigma \hookrightarrow
  S^2\setminus P$. Then the \emph{inverse image} $f^{-1}(\Sigma)$ is the
  pull-back in the diagram
  \[
  \begin{tikzpicture}
    \matrix[row sep=0.8cm,column sep=1cm] {
      \node (Sigmai) {$f^{-1}(\Sigma)$}; &
        \node (Sigma) {$\Sigma$}; \\
      \node (S2i) {$S^2\setminus f^{-1}(P)$}; &
        \node (S2) {$S^2\setminus P$.}; \\
    };
    \draw[->] (Sigma) to node[auto=left,cdlabel] {i} (S2);
    \draw[->] (S2i) to node[auto=left,cdlabel] {f} (S2);
    \draw[->] (Sigmai) to node[auto=left,cdlabel] {\pi} (Sigma);
    \draw[->] (Sigmai) to (S2i);
  \end{tikzpicture}
  \]
  There are two natural maps:
  \begin{itemize}
  \item An embedding $f^{-1}(\Sigma) \hookrightarrow S^2\setminus
    f^{-1}(P) \subset S^2\setminus P$.
  \item A covering map $\pi\co f^{-1}(\Sigma) \to \Sigma$. Since it is a
    covering map, $f^{-1}(\Sigma)$ inherits all the structure of $\Sigma$.
  \end{itemize}
\end{definition}

If $\Sigma$ is a surface spine for $S^2 \setminus P$, then the
embedding $f^{-1}(\Sigma) \to S^2 \setminus P$ is isotopic to
a topological embedding $\phi \co f^{-1}(\Sigma) \to \Sigma$. We
therefore have a virtual endomorphism of~$\Sigma$ and can
apply the technology of Section~\ref{sec:dynamics}.

\begin{citethm}\label{thm:rational-surfaces-embed}
  Let $f\co (S^2,P)\righttoleftarrow$ be a branched self-cover with
  a branch point in each cycle. Then $f$ is equivalent to a rational
  map iff there is a
  conformal surface spine~$\Sigma$ for $S^2\setminus P$ so that
  in the corresponding virtual endomorphism, $\phi$ is homotopic to a
  strict conformal embedding $f^{-1}(\Sigma) \hookrightarrow \Sigma$.
\end{citethm}

The above characterization of rational maps appears to have been
folklore in the community for some time.

\begin{proof}[Proof sketch]
  If $f\co(\CC\PP^1,P)\righttoleftarrow$ is a rational map, then we
  may take $\Sigma$ to be a neighborhood of the Julia set
  $J(f)$ (filling in all disks that do not contain critical
  points). Since the dynamics is super-attracting to~$P$ on the
  complement of~$J(f)$, the inverse image $f^{-1}(\Sigma)$ is a
  smaller neighborhood of $J(f)$ than $\Sigma$, and by taking $\Sigma$
  to be suitably
  nice (e.g., using level sets of the potential function) we can
  arrange for $f^{-1}(\Sigma) \subset \Sigma$ to be a strict conformal
  embedding.

  The converse direction is a special case of \cite[Theorem
  5.2]{CPT14:Renorm} or \cite[Theorem
  7.1]{Wang14:DecompositionHerman}. The technique, quasi-conformal
  surgery,
  goes back to Douady and Hubbard \cite{DH85:DynamicsPolyLike}.
\end{proof}

\begin{remark}\label{rem:embed-stronger}
  In fact, Theorem~\ref{thm:rational-surfaces-embed} is true if the
  assumption that the conformal embedding is strict is dropped. With
  that modification, it also extends to more general branched
  self-covers. This strengthening is not needed for the main theorem.
\end{remark}

\subsection{Characterizing rational maps using graphs}
\label{sec:rational-graphs}

In the setting of Theorem~\ref{thm:embedding-graph-surface}, we can
consider the asymptotic stretch factor of the virtual endomorphism.
Theorem~\ref{thm:rational-surfaces-embed} and
Propositions~\ref{prop:asf-bound-surf} and
\ref{prop:asf-to-sf} then tell us the
following.

\begin{proposition}
  Let $f \co (S^2,P)\righttoleftarrow$ be a branched self-cover with a
  branch point in each cycle in~$P$ and let $\pi,\phi\co \Sigma_1
  \rightrightarrows \Sigma$ be the
  corresponding virtual endomorphism. Then $f$ is equivalent to a
  rational map iff $\ASF[\phi] < 1$.
\end{proposition}

By the equality of the asymptotic stretch factors for graphs and for
surfaces (Proposition~\ref{prop:asf-graph-surface}), we immediately
get the following more precise version of Theorem~\ref{thm:detect-rational}.

\begin{theorem}\label{thm:detect-rational-2}
  Let $[f] \co (S^2,P) \righttoleftarrow$ be a branched self-cover of
  the sphere relative to a finite number of points $P \subset
  S^2$, with a branch point in each cycle
  in~$P$. For $\Gamma$ any spine of
  $S^2\setminus P$, let $\pi_\Gamma,\phi_\Gamma\co \Gamma_1
  \rightrightarrows \Gamma$ be the
  corresponding virtual endomorphism. Then
  the following conditions are equivalent.
  \begin{enumerate}
  \item\label{item:cond-rational} The branched self-cover $f$ is
    equivalent to a rational map.
  \item\label{item:cond-asf} For any spine~$\Gamma$ for $S^2\setminus
    P$, we have $\ASF[\phi_\Gamma] < 1$.
  \item\label{item:cond-exist-Gamma} There is some elastic graph
    spine~$\Gamma$ for $S^2\setminus
    P$ and some integer $n > 0$ so that $\Emb[\phi_{\Gamma,n}] < 1$.
  \item\label{item:cond-all-Gamma} For every elastic graph
    spine~$\Gamma$ for $S^2\setminus
    P$ and every sufficiently large $n$, we have $\Emb[\phi_{\Gamma,n}] < 1$.
  \end{enumerate}
\end{theorem}

\begin{proof}[Proof sketch]
The proof has essentially already been given earlier in the paper.
The diagram below summarizes the chain of implications between
(\ref{item:cond-rational}), (\ref{item:cond-asf}),
and~(\ref{item:cond-exist-Gamma}). Here, $\Sigma$ is a (conformal)
surface spine for $S^2 \setminus P$ and $\phi_\Sigma$ is the
corresponding surface virtual endomorphism. There are various side
conditions for which you should see the referenced theorems and
propositions.

\begin{equation*}
\begin{tikzpicture}
  \matrix[row sep=6ex, column sep=2em] {
    \node (rational) {$f$ rational}; \\
    \node (emb-surf) {\minibox[c]{%
        $\exists \Sigma$: $f^{-1}(\Sigma)$ strictly\\
        conformally embeds in $\Sigma$}}; &
      \node (emb-graph) {$\exists\Gamma, n$: $\Emb[\phi_{\Gamma,n}] < 1$}; \\
    \node (SF-surf) {$\exists \Sigma$: $\SF[\phi_\Sigma] < 1$}; &
      \node (SF-graph) {$\exists\Gamma, n$: $\SF[\phi_{\Gamma,n}] < 1$}; \\
    \node (ASF-surf) {$\ASF[\phi_\Sigma] < 1$}; &
      \node (ASF-graph) {$\ASF[\phi_\Gamma] < 1$}; \\
  };
  \draw[biimplication] (rational) to
    node[left]{\footnotesize{Thm.\ \ref{thm:rational-surfaces-embed}}}
    (emb-surf);
  \draw[biimplication] (emb-surf) to
    node[left]{\footnotesize{Thm.\ \ref{thm:el-embedding}}}
    (SF-surf);
  \draw[implication,bend right] (SF-surf) to
    node[left]{\footnotesize{Prop.\ \ref{prop:asf-bound-surf}}}
    (ASF-surf);
  \draw[implication,bend right] (ASF-surf) to
    node[right]{\footnotesize{\minibox[l]{Prop.\ \ref{prop:asf-to-sf},\\
        See below}}}
    (SF-surf);
  \draw[biimplication] (ASF-surf) to
    node[auto=right]{\footnotesize{Prop.\ \ref{prop:asf-graph-surface}}}
    (ASF-graph);
  \draw[implication,bend right] (SF-graph) to
    node[left]{\footnotesize{Prop.\ \ref{prop:asf-bound-graph}}}
    (ASF-graph);
  \draw[implication,bend right] (ASF-graph) to
    node[right]{\footnotesize{Prop.\ \ref{prop:asf-to-sf}}}
    (SF-graph);
  \draw[biimplication] (emb-graph) to
    node[right]{\footnotesize{Thm.\ \ref{thm:sf-emb-energy}}}
    (SF-graph);
\end{tikzpicture}
\end{equation*}
Lemma~\ref{lem:asf-indep} says that the asymptotic stretch factor
$\ASF$ is independent of the conformal structure on~$\Sigma$ or
elastic structure on~$\Gamma$, giving the equivalence to
statement~(\ref{item:cond-all-Gamma}).

The implication
\[
\ASF[\phi_\Sigma] < 1 \Rightarrow
\exists\Sigma\co \SF[\phi_\Sigma] < 1
\]
is not direct. Proposition~\ref{prop:asf-to-sf} gives the implication
\[
\ASF[\phi_\Sigma] < 1 \Rightarrow \exists\Sigma,n\co \SF[\phi_{\Sigma,n}] < 1,
\]
which by Theorem~\ref{thm:rational-surfaces-embed} implies that
$f^{\circ n}$
is equivalent to a rational map. To then conclude that $f$ is a
rational map and that we don't need to iterate to find a surface
spine~$\Sigma$ with $\SF[\phi_\Sigma] < 1$,
we need to know a little bit more about the
geometry. For instance, it suffices to know that, except in Lattés
examples, some power of the backwards iteration map on Teichmüller
space is contracting. This is the only place that this argument relies
on the original characterization of rational maps~\cite{DH93:ThurstonChar}.
\end{proof}

The virtual endomorphism $\phi \co f^{-1}(\Gamma) \to \Gamma$ is going in
essentially the reverse direction to~$f$. While $f$ is expanding on
the Julia set, $\phi$ is contracting (in an appropriate
sense). Therefore, for $f$ a branched self-cover with a branch point
in each cycle and $\phi$ the associated
virtual endomorphism of a spine, we define
\begin{equation}
  \ASF[f] \coloneqq \ASF[\phi]^{-1}.
  \label{eq:def-asf-ratl}
\end{equation}
This quantity is greater than one when $f$ is rational.

\subsection{Polynomials}
\label{sec:polynomials}

One important special case is that of polynomials, where
Theorem~\ref{thm:detect-rational-2} recovers known results, which we
briefly review.

\begin{definition}
  A \emph{topological polynomial} is a branched self-cover $[f]\co
  (S^2,P)\righttoleftarrow$ of degree~$d$ with one fixed point $\infty
  \in P$ which is branched of degree~$d$. If such a map is equivalent
  to a rational map, the rational map is a polynomial.
\end{definition}

\begin{definition}
  For a post-critically finite polynomial~$p$ with post-critical set $P
  \cup \{\infty\}\subset\CC\PP^1$, the \emph{filled Julia set} $K(p)$
  is the set of points bounded under iteration. It is the
  union of the Julia set~$J(p)$ and the attracting basin around each
  finite \emph{Fatou point}, points that eventually maps to a periodic cycle
  with a branch point.
  The \emph{Hubbard tree}~$T_P\subset K(p)$ is
  essentially the spanning tree of~$P$ within $K(p)$. When $T_P$
  intersects a Fatou component, it is required to be a union of rays
  in the canonical coordinates. See
  \cite{DH85:DynamicsPolyLike}.

  Now define the \emph{Hubbard graph} $G_P\subset J(p)$ to be the
  union of
  \begin{itemize}
  \item  the boundary of the Fatou component containing each Fatou
    point in~$P$,
  \item the Julia points of~$P$, and
  \item the minimal tree connecting the above points in $K(p)$, taking
    unions of rays in coordinates of Fatou components as before.
  \end{itemize}
\end{definition}

$T_P$ is forward invariant under the polynomial~$p$. As such, we can
consider dynamics
on it.
To make it a complete invariant of the map, we can decorate $T_P$
with some additional angle data at Fatou points \cite[Section
6.1.2]{DH84:OrsayI}; this is close to what
is needed to reconstruct the graph~$G_P$ from the map on~$T_P$.
A \emph{abstract Hubbard tree} is a tree with an endomorphism, marked
with this additional data.

\begin{definition}
  An endomorphism $f$ of a Hubbard tree is \emph{Julia expanding} if it
  admits a metric for which $f$ does not decrease distances and
  strictly increases the
  distance between any two Julia points.
\end{definition}

The Hubbard graph $G_P$, on the other hand, is not in general forward
invariant, since near a point in $f^{-1}(P)\setminus P$ the
graph~$G_P$ may go
through a Fatou component. It does
have a natural virtual endomorphism, which can be
reconstructed from the abstract Hubbard tree.

\begin{lemma}
  An endomorphism of an abstract Hubbard tree is Julia expanding iff
  the corresponding virtual endomorphism of the abstract Hubbard graph
  has embedding energy less than~$1$.
\end{lemma}

\begin{citethm}[Poirier \cite{Poirier09:CritPortraits,Poirier10:HubbardTrees}]
  \label{thm:realize-hubbard}
  An abstract Hubbard tree $T_P$ is realizable as the Hubbard tree of
  a polynomial iff it is Julia expanding.
\end{citethm}

In the presence of expanding dynamics, we can consider the entropy.

\begin{definition}
  The \emph{core entropy} $h(f)$ of a post-critically finite
  polynomial~$f$ is
  the topological entropy of the map on the Hubbard tree~$T_P$.
\end{definition}

\begin{proposition}
  Let $f$ be a post-critically finite quadratic polynomial. If
  $f$ is dendritic (i.e., each critical point is strictly
  pre-periodic), then
  \[\log \ASF[f] = h(f).\]
\end{proposition}
The statement uses an extension of the theory of stretch factors,
etc., to the non-hyperbolic
case by using marked orbifold points, as mentioned in
Problem~\ref{prob:non-hyp} below.

\begin{example}\label{examp:z2+i}
  Consider the polynomial $f(z) = z^2 + i$. The critical portrait
  is
  \[
  \mathcenter{\begin{tikzpicture}[node distance=1.75cm]
      \node (A) at (0,0) {$0$};
      \node (B)[right of=A] {$i$};
      \node (C)[right of=B] {$-1+i$};
      \node (D)[right of=C] {$-i$};
      \node (E)[right of=D] {$\infty$};
      \draw[|->] (A) to node[above,cdlabel]{(2)} (B);
      \draw[|->] (B) to (C);
      \draw[|->, bend left=25] ($(C.east)+(0,0.75ex)$) to ($(D.west)+(0,0.75ex)$);
      \draw[|->, bend left=25] ($(D.west)-(0,0.75ex)$) to ($(C.east)-(0,0.75ex)$);
      \draw[|->, loop above, min distance=0.8cm,out=110,in=70] (E) to
          node[above,cdlabel]{(2)} (E);
    \end{tikzpicture}}
  \]
\begin{figure}
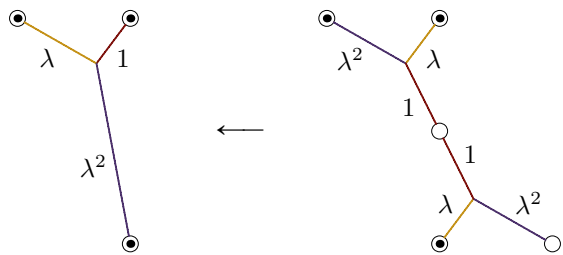

  \[
  \mfigb{graphs-50} \quad\quad\longleftarrow\quad\,\, \mfigb{graphs-51}
  \]
  \caption{Spines for the dendritic polynomial $z^2 + i$}
  \label{fig:z2+i}
\end{figure}
  A spine $\Gamma$ and its inverse image $\wt\Gamma$ for this
  polynomial are shown in
  Figure~\ref{fig:z2+i}. In the extended theory with marked orbifold
  points, the small loops in Figure~\ref{fig:z2+i} should be thought
  of as being of length~$0$, and ignored in the computation of the
  stretch factor. If we do this and set $\lambda\approx 1.52$ to be the
  positive root of $\lambda^3 = 2 + \lambda$, then the stretch factor
  of the map $\phi \co \wt\Gamma \to \Gamma$ is visibly
  $1/\lambda$. This polynomial equation is exactly the same equation
  one solves to find the core entropy.
\end{example}

  If $f$ is hyperbolic (i.e.,
  each critical point ends up in a cycle with at least one critical point)
  and the critical cycles have degree $d_i$ and length $n_i$ (for $i=1,\dots,k$)
  then
  \[\log \ASF[f] \ge \max_{1 \le i \le k} \frac{\log d_i}{n_i}.\]
  This last inequality is often an equality and often gives answers
  larger than $h(f)$.

The fact that $f$ is realized as a polynomial iff the
entropy on a topological Hubbard graph is positive is essentially a
corollary of Poirier's Theorem~\ref{thm:realize-hubbard}.

\subsection{Comparison to annular obstruction}
\label{sec:annular}

As mentioned in the introduction, W.~Thurston in 1982 gave a different
characterization of rational maps among topological branched
self-covers, as explained by Douady and Hubbard. Instead of finding an
object that guarantees that the map is rational, he finds an
obstruction that guarantees it is not rational. We can write down
that obstruction in our language as follows.

\begin{definition}
  An \emph{annular} elastic graph is one that consists only of circle
  components. Let $A$ be an annular elastic graph embedded in an
  orbifold~$\Sigma$. Then the \emph{join} of~$A$ is the elastic
  graph $\Join(A)$ obtained by
  \begin{itemize}
  \item deleting all components of~$A$ that bound a disk with at most
    one orbifold point, and then
  \item merging all parallel components of~$A$ by harmonically adding
    their elastic constants.
  \end{itemize}
  \emph{Harmonically adding} the elastic constants means taking two parallel
  components $a_1$ and~$a_2$ and replacing them with a new component~$a_3$
  with
  \[
  \frac{1}{\alpha(a_3)} = \frac{1}{\alpha(a_1)} + \frac{1}{\alpha(a_2)}.
  \]
  This corresponds to stitching together two side-by-side conformal
  annuli in the most efficient way possible.
\end{definition}

\begin{citethm}[W.~Thurston, Douady-Hubbard \cite{DH93:ThurstonChar}]
  \label{thm:thurston-obstruction}
  Let $f\co (S^2,P)\righttoleftarrow$ be a topological branched
  self-cover that is not a Lattés example. Let $S^2_f$ be the orbifold
  of~$f$. Then $f$ is
  equivalent to a rational map iff there is no annular elastic
  graph~$A$ in $S^2_f$ and map $\phi \co A \to \Join(f^{-1}(A))$
  compatible with the maps to~$S^2\setminus P$ so that $\Emb(\phi) \le 1$.
\end{citethm}

\begin{remark}
  The usual formulation of
  Theorem~\ref{thm:thurston-obstruction} refers to the
  maximum eigenvalue of a matrix constructed out of
  $A$ considered as a multi-curve (with no extra structure). The
  above formulation is easily seen to be equivalent.
\end{remark}

Where our Theorem~\ref{thm:detect-rational-2} looks for a strictly
loosening map
\begin{align*}
f^{-1}(\Gamma) &\rightarrow \Gamma,
\intertext{the older Theorem~\ref{thm:thurston-obstruction} looks for a (not
necessarily strictly) loosening map the other direction}
\Join(f^{-1}(A)) &\leftarrow A.
\end{align*}
There is an easy argument that both conditions cannot simultaneously
hold (as implied by the theorems). If both loosening maps
existed, scale
$A$ so that $\Emb[A \to \Gamma] = 1$, and consider the
square of loosening maps
\begin{equation}
\begin{tikzpicture}[baseline]
  \matrix[row sep=0.8cm, column sep=1cm]{
    \node (Ai) {$\Join(f^{-1}(A))$}; &
      \node (A) {$A$}; \\
    \node (Gi) {$f^{-1}(\Gamma)$}; &
      \node (G) {$\Gamma.$}; \\
  };
  \draw[->] (A) to (Ai);
  \draw[->] (Ai) to (Gi);
  \draw[->] (Gi) to (G);
  \draw[->] (A) to (G);
\end{tikzpicture}\label{eq:obstruction-cd}
\end{equation}
The map $f^{-1}(\Gamma) \hookrightarrow \Gamma$ is strictly
loosening but the map $A \hookrightarrow \Gamma$ is not, a
contradiction.

It is also worth noting that
Theorem~\ref{thm:detect-rational-2} is technically easier than
the existing proof of Theorem~\ref{thm:thurston-obstruction}. The key
analytic point is contained in
Theorem~\ref{thm:rational-surfaces-embed}, which is relatively soft.

\subsection{Comparison to domination of weighted arc diagrams}
\label{sec:dynam-teich}

Suppose we are given a virtual endomorphism of surfaces $\pi,\phi\co
\Sigma_1 \rightrightarrows \Sigma_0$,
with $\pi$ a covering map and $\phi$ a topological embedding.
Then the associated \emph{dynamical Teichm\"uller space} is
\begin{equation*}
\Teich(\pi,\phi) \coloneqq \bigl\{\,S \in \Teich(\Sigma_0)\bigm|
   \text{$\pi^*S$ conformally embeds in $S$ in the homotopy class $[\phi]$}\,
   \bigr\}.
\end{equation*}
Here, $\Teich(\Sigma_0)$ is the finite-dimensional space of conformal
structures on the interior of~$\Sigma_0$.  (The conformal structures
are allowed to have removable singularities or not; that is, the
corresponding hyperbolic structures are allowed to have parabolic or
hyperbolic monodromy
around the boundary components.)  Then
Theorem~\ref{thm:rational-surfaces-embed} and the strengthening in
Remark~\ref{rem:embed-stronger} say that the virtual endomorphism
comes from a rational map iff $\Teich(\pi,\phi)$ is non-empty.

Jeremy Kahn has studied the topology of $\Teich(\pi,\phi)$ in the
context of renormalization of polynomials, and has results in terms of
\emph{weighted arc
  diagrams} \cite{Kahn06:BoundsI}.
A weighted arc diagram is a weighted collection~$X$ of arcs with
endpoints on boundary
components of a surface~$\Sigma$ with boundary. If $X$ is filling (the
complementary components are disks), it can be thought of as dual to
an elastic graph~$\Gamma(X)$ embedded
in $\Sigma$, with one
vertex per component of $\Sigma\setminus X$ and one edge crossing
each arc in~$X$. The elastic length of an edge of $\Gamma(X)$ is
equal to the weight on the arc it crosses.

Now suppose that you have two surfaces, $\Sigma_1 \subset
\Sigma_0$. There if $X$ is a weighted arc diagram on~$\Sigma_1$ and
$Y$ is a weighted arc diagram on~$\Sigma_0$, there is a notion of $X$
\emph{dominating} $Y$, written $X \multimap Y$. We do not repeat the
definition here (see \cite[Section
3.6]{Kahn06:BoundsI}), but, in the filling case, $X$ dominates $Y$ iff
there is a map $\phi
\co \Gamma(X) \to \Gamma(Y)$
with $\Emb(\phi) \le 1$ that is compatible up to homotopy with the
inclusion $\Sigma_1 \subset \Sigma_0$.

Kahn used this notion of domination to control certain infinitely
renormalizable polynomials (studying how conformal structures can
degenerate). In terms of Teichm\"uller spaces, we have the following:

\begin{proposition}[Kahn, personal communication]
  If $\Teich(\pi,\phi)$ is not compact, then there is a weighted arc
  diagram $X$ on $\Sigma_0$ so that $\pi^*X \multimap X$.
\end{proposition}

\begin{conjecture}[Kahn, personal communication]\label{conj:wad-dom-noncompact}
  If there is a weighted arc diagram $X$ on $\Sigma_0$ so that $\pi^*X
  \multimap X$, then $\Teich(\pi,\phi)$ is not compact.
\end{conjecture}

We can prove a weaker version of Conjecture~\ref{conj:wad-dom-noncompact}:
\begin{corollary}\label{cor:strict-dom-noncompact}
  If there is a weighted arc diagram $X$ on $\Sigma_0$ so that
  $\pi^* X$ strictly dominates $X$ and
  $X$ fills $\Sigma_0$, then $\Teich(\pi,\phi)$ is not compact.
\end{corollary}
\begin{proof}
  The conditions are equivalent to saying that $\Gamma(X)$ is a spine
  for $\Sigma_0$ and $\Emb[\pi^* \Gamma(X)\to \Gamma(X)] < 1$. Then
  Theorem~\ref{thm:embedding-graph-surface} says that
  $N_t\Gamma(X)\in \Teich(\pi,\phi)$ for $t$ sufficiently small; this
  gives an explicit sequence exhibiting the non-compactness of
  $\Teich(\pi,\phi)$.
\end{proof}

\begin{corollary}
  If the virtual endomorphism $\pi,\phi\co \Sigma_1 \to \Sigma$ comes
  from a topological
  branched self-cover $f\co (S^2,P)\righttoleftarrow$ with at least
  one branch point in each cycle, then $\Teich(\pi_n,\phi_n)$ is not
  compact for $n$ sufficiently large.
\end{corollary}

\begin{proof}
  Theorem~\ref{thm:detect-rational-2} says that for $n$ large enough,
  there is an elastic graph~$\Gamma$ so that $\Emb[f^{-n}(\Gamma)\to
  \Gamma] < 1$, which by Corollary~\ref{cor:strict-dom-noncompact}
  implies $\Teich(\pi_n,\phi_n)$ is not compact.
\end{proof}

In fact, the techniques can be strengthened to show that as long as
$\Teich(\pi,\phi)$ is more than a single point, then
$\Teich(\pi_n,\phi_n)$ is eventually non-compact.

In addition, Kahn has communicated examples where
$\Teich(\pi,\phi)$ is compact (but more than one point), thus showing
that iteration is
necessary in Theorem~\ref{thm:detect-rational-2}.

\subsection{Comparison to Lipschitz expansion}
\label{sec:branched-lipschitz}

In Section~\ref{sec:asf-lipschitz} we remarked that there are two
analytic quantities associated to a virtual endomorphism of a graph: the
asymptotic stretch factor $\ASF$ and the asymptotic Lipschitz constant
$\ALip$. In the setting of Theorem~\ref{thm:detect-rational-2}, what
are the virtual endomorphisms for which $\ALip[f^{-1}(\Gamma) \to \Gamma] < 1$? Since
$\ALip[\phi] \le \ASF[\phi]$
(Equation~\eqref{eq:ASF-ALip}), this is a weaker
condition. In fact, $\ALip[\phi] < 1$ can be taken as a definition of
when
associated virtual endomorphism of $\pi_1(S^2\setminus P)$ is
contracting in the sense of Nekrashevych
\cite{Nekrashevych05:SelfSimilar,nekrashevych11:IMG,Nekrashevych13:CombModel}.

Bonk and Meyer \cite{BM10:Expanding} and Haïssinsky and Pilgrim
\cite{HP12:AlgExpanding} studied
essentially this condition, in the opposite situation to the present paper:
they considered the case when there are
no branch points in any cycle in~$P$. (For rational maps, this means
that the Julia set is the whole sphere.) In our language, they showed
more-or-less that for a branched self-cover~$f$ with no branch points
in periodic cycles with associated virtual
endomorphism~$\phi$, the stretch factor
 $\ALip[\phi] < 1$ iff $f$ has a
topological representative as a map on $S^2$ that is uniformly
expanding. Any rational map has such a uniformly expanding
representative, and many other maps do as
well.

\subsection{Extensions and questions}
\label{sec:rational-extensions}

The first obvious question is to extend
Theorem~\ref{thm:detect-rational-2} away from the restricted setting.

\begin{problem}\label{prob:non-hyp}
  Extend Theorem~\ref{thm:detect-rational-2} to general topological
  branched self-covers, dropping the condition that there be a
  critical point in each cycle.
\end{problem}

It appears that the theory will extend without much problem to the
much larger family of maps which have at least one critical point in
one cycle in~$P$. (For rational maps, these are maps where the Julia
set is not the entire sphere.) This requires extending the theory for
graphs and surfaces to allow dealing with marked orbifold points, as
briefly
mentioned in Section~\ref{sec:extensions}.

For the remaining case, when there are no branch points in any cycles
in~$P$ and the Julia set is the whole sphere if the map is rational,
there is no clear analogue of the conformal
characterization in Theorem~\ref{thm:rational-surfaces-embed} and so it
is not clear what the statement would be.

\begin{question}\label{quest:n=1}
  In Theorem~\ref{thm:detect-rational-2}, does $n=1$ suffice? That is,
  for $f$ as given there, is there always an elastic graph~$\Gamma$ so
  that $\Emb[f^{-1}(\Gamma) \to \Gamma] < 1$?
\end{question}

The answer to Question~\ref{quest:n=1} is ``yes'' for
polynomials and for many examples, and
Theorem~\ref{thm:rational-surfaces-embed} is a corresponding statement for
conformal surfaces. Nevertheless, the answer
to Question~\ref{quest:n=1} is ``no'' in
general. See
Section~\ref{sec:dynam-teich}.

\begin{problem}\label{prob:min-n}
  In Theorem~\ref{thm:detect-rational-2}, how does the minimal
  necessary $n$ grow, as a function of the degree and the size of~$P$?
\end{problem}

\begin{problem}\label{prob:alg-complexity}
  Study the algorithmic complexity of finding a strictly loosening map
  $\phi \co f^{-n}(\Gamma) \to \Gamma$ when $f$ is
  equivalent to a rational map. In particular, does this give a
  polynomial-sized certificate that $f$ is rational?
\end{problem}

Presumably a solution to Problem~\ref{prob:alg-complexity} would first
require a solution to Problem~\ref{prob:min-n}.

\begin{problem}
  Study $\ASF[f]$ as $f$ varies over rational maps. How does it
  compare to known invariants?
\end{problem}

Question~\ref{quest:asf-algebraic} is one starting question.
The connection between $\ASF$ and the entropy on the Hubbard graph
described in Section~\ref{sec:polynomials} is relevant, but in general
there is no
reasonable forward dynamics on~$\Gamma$.

Finally, Theorems~\ref{thm:detect-rational-2}
and~\ref{thm:thurston-obstruction} give opposite combinatorial
conditions. These should presumably be related.

\begin{problem}\label{prob:comb-annulus}
  Prove combinatorially that, for any virtual endomorphism $\phi \co
  \wt\Gamma \to \Gamma$ of a graph, either $\ASF[\phi] < 1$, or there
  is an annular obstruction in the sense of
  Theorem~\ref{thm:thurston-obstruction}.
\end{problem}

A solution to Problem~\ref{prob:comb-annulus} would give an
alternate proof of the main Theorem~\ref{thm:detect-rational-2}, using
W.~Thurston's
Theorem~\ref{thm:thurston-obstruction}.

Another family of problems comes from applying
Theorem~\ref{thm:detect-rational-2}: Prove that certain maps
\emph{are} rational by finding explicit maps $\phi \co f^{-1}(\Gamma)
\to \Gamma$ with $\Emb(\phi) < 1$.  We give one sample problem.

\begin{problem}
  Rees and Tan Lei \cite{Tan92:Matings}, and later Shishikura
  \cite{Shishikura00:ReesMatings} showed that two quadratic
  polynomials are topologically mateable iff they do not lie in
  conjugate limbs of the Mandelbrot set. Reprove this with the
  techniques of this paper. (See
  Examples~\ref{examp:mate-rabbit-basilica} and
  \ref{examp:mate-basilica-basilica} for one approach.) More
  generally, find a criterion for cubic or
  higher-degree polynomials to be mateable. This is likely to be more
  difficult, given the cubic example of Shishikura and Tan Lei without
  Levy cycles \cite{ST00:CubicMating}.
\end{problem}

\begin{warning}
  A topological mating is not the naive (formal) mating obtained by
  gluing two polynomials at infinity. See the survey by Buff et
  al.\ \cite{BEKMPRL12:PolyMatings}.
\end{warning}


\bibliographystyle{hamsalpha}
\bibliography{summary}

\end{document}